# Properties of distance spaces with power triangle inequalities

Daniel J. Greenhoe

*Abstract*: *Metric space*s provide a framework for analysis and have several very useful properties. Many of these properties follow in part from the *triangle inequality*. However, there are several applications in which the triangle inequality does not hold but in which we may still like to perform analysis. This paper investigates what happens if the triangle inequality is removed all together, leaving what is called a *distance space*, and also what happens if the triangle inequality is replaced with a much more general two parameter relation, which is herein called the "*power triangle inequality*". The power triangle inequality represents an uncountably large class of inequalities, and includes the *triangle inequality*, *relaxed triangle inequality*, and *inframetric inequality* as special cases. The power triangle inequality is defined in terms of a function that is herein called the *power triangle function*. The power triangle function is itself a *power mean*, and as such is *continuous* and *monotone* with respect to its exponential parameter, and also includes the operations of *maximum*, *minimum*, *mean square*, *arithmetic mean*, *geometric mean*, and *harmonic mean* as special cases.

# Contents

# 1    Introduction and summary

*Metric space*s provide a framework for analysis and have several very useful properties. Many of these properties follow in part from the *triangle inequality*. However, there are several applications[1] in which the triangle inequality does not hold but in which we would still like to perform analysis. So the questions that natually follow are:

Q1.   What happens if we remove the *triangle inequality* all together?
Q2.   What happens if we replace the *triangle inequality* with a generalized relation?

A *distance space* is a *metric space* without the *triangle inequality* constraint. Section 3 introduces *distance space*s and demonstrates that some properties commonly associated with *metric space*s also hold in any *distance space*:

D1.   ∅ and $X$ are *open*                                        (Theorem 3.7 page 7)
D2.   the intersection of a finite number of open sets is *open*   (Theorem 3.7 page 7)
D3.   the union of an arbitrary number of open sets is *open*      (Theorem 3.7 page 7)
D4.   every Cauchy sequence is *bounded*                           (Proposition 3.14 page 9)
D5.   any subsequence of a *Cauchy* sequence is also *Cauchy*      (Proposition 3.15 page 10)
D6.   the *Cantor Intersection Theorem* holds                      (Theorem 3.18 page 11)

---

[1]references for applications in which the *triangle inequality* may not hold: 📖 [Maligranda and Orlicz(1987)] page 54 ⟨"pseudonorm"⟩, 📖 [Lin(1998)] ⟨"similarity measures", Table 6⟩, 📖 [Veltkamp and Hagedoorn(2000)] ⟨"shape similarity measures"⟩, 📖 [Veltkamp(2001)] ⟨"shape matching"⟩, 📖 [Costa et al.(2004)Costa, Castro, Rowstron, and Key] ⟨"network distance estimation"⟩, 📖 [Burstein et al.(2005)Burstein, Ulitsky, Tuller, and Chor] page 287 ⟨distance matrices for "genome phylogenies"⟩, ✏ [Jiménez and Yukich(2006)] page 224 ⟨"statistical distances"⟩, 📖 [Szirmai(2007)] page 388 ⟨"geodesic ball"⟩, 📖 [Crammer et al.(2007)Crammer, Kearns, and Wortman] page 326 ⟨"decision-theoretic learning"⟩, 📖 [Vitányi(2011)] page 2455 ⟨"information distance"⟩,

The following five properties (M1–M5) *do* hold in any *metric space*. However, the examples from Section 3 listed below demonstrate that the five properties do *not* hold in all *distance spaces*:

| | | | |
|---|---|---|---|
| M1. | the *metric function* is *continuous* | fails to hold in | Example 3.21–Example 3.23 |
| M2. | *open ball*s are *open* | fails to hold in | Example 3.21 and Example 3.22 |
| M3. | the *open ball*s form a *base* for a topology | fails to hold in | Example 3.21 and Example 3.22 |
| M4. | the limits of *convergent sequences* are *unique* | fails to hold in | Example 3.21 |
| M5. | *convergent* sequences are *Cauchy* | fails to hold in | Example 3.22 |

Hence, Section 3 answers question Q1.

Section 4 begins to answer question Q2 by first introducing a new function, called the *power triangle function* in a *distance space* $(X, \mathsf{d})$, as $\quad \tau(p,\sigma;x,y,z;\mathsf{d}) \triangleq 2\sigma\left[\frac{1}{2}\mathsf{d}^p(x,z) + \frac{1}{2}\mathsf{d}^p(z,y)\right]^{\frac{1}{p}} \quad$ for some $(p,\sigma) \in \mathbb{R}^* \times \mathbb{R}$.[2] Section 4 then goes on to use this function to define a new relation, called the *power triangle inequality* in $(X, \mathsf{d})$, and defined as $\quad \bigcirc(p,\sigma;\mathsf{d}) \triangleq \left\{(x,y,z) \in X^3 \mid \mathsf{d}(x,y) \le \tau(p,\sigma;x,y,z;\mathsf{d})\right\}$.

The *power triangle inequality* is a generalized form of the *triangle inequality* in the sense that the two inequalities coincide at $(p,\sigma) = (1,1)$. Other special values include $(1,\sigma)$ yielding the *relaxed triangle inequality* (and its associated *near metric space*) and $(\infty,\sigma)$ yielding the *$\sigma$-inframetric inequality* (and its associated *$\sigma$-inframetric space*). Collectively, a distance space with a power triangle inequality is herein called a *power distance space* and denoted $(X, \mathsf{d}, p, \sigma)$.[3]

The *power triangle function*, at $\sigma = \frac{1}{2}$, is a special case of the *power mean* with $N = 2$ and $\lambda_1 = \lambda_2 = \frac{1}{2}$. *Power mean*s have the elegant properties of being *continuous* and *monontone* with respect to a free parameter $p$. From this it is easy to show that the *power triangle function* is also *continuous* and *monontone* with respect to both $p$ and $\sigma$. Special values of $p$ yield operators coinciding with *maximum*, *minimum*, *mean square*, *arithmetic mean*, *geometric mean*, and *harmonic mean*. *Power mean*s are briefly described in APPENDIX B.2.[4]

Section 4.2 investigates the properties of *power distance spaces*. In particular, it shows for what values of $(p,\sigma)$ the properties M1–M5 hold. Here is a summary of the results in a *power distance space* $(X, \mathsf{d}, p, \sigma)$, for all $x, y, z \in X$:

| | | |
|---|---|---|
| (M1) | holds for any $(p,\sigma) \in (\mathbb{R}^*\setminus\{0\}) \times \mathbb{R}^+$ such that $2\sigma = 2^{\frac{1}{p}}$ | (Theorem 4.18 page 23) |
| (M2) | holds for any $(p,\sigma) \in (\mathbb{R}^*\setminus\{0\}) \times \mathbb{R}^+$ such that $2\sigma \le 2^{\frac{1}{p}}$ | (Corollary 4.14 page 21) |
| (M3) | holds for any $(p,\sigma) \in (\mathbb{R}^*\setminus\{0\}) \times \mathbb{R}^+$ such that $2\sigma \le 2^{\frac{1}{p}}$ | (Corollary 4.12 page 20) |
| (M4) | holds for any $(p,\sigma) \in \mathbb{R}^* \times \mathbb{R}^+$ | (Theorem 4.19 page 23) |
| (M5) | holds for any $(p,\sigma) \in \mathbb{R}^* \times \mathbb{R}^+$ | (Theorem 4.16 page 21) |

APPENDIX A briefly introduces *topological space*s. The *open ball*s of any *metric space* form a *base* for a *topology*. This is largely due to the fact that in a metric space, open balls are *open*. Because of this, in metric spaces it is convenient to use topological structure to define and exploit analytic concepts such as *continuity*, *convergence*, *closed set*s, *closure*, *interior*, and *accumulation point*. For example, in a metric space, the traditional definition of defining continuity using open balls and the topological

---

[2] where $\mathbb{R}^*$ is the *set of extended real numbers* and $\mathbb{R}^+$ is the *set of positive real numbers* (Definition 2.1 page 4)

[3] *power triangle inequality*: Definition 4.4 page 16; *power distance space*: Definition 4.3 page 16; examples of *power distance space*: Definition 4.5 page 16;

[4] *power triangle function*: Definition 4.1 (page 15); *power mean*: Definition B.6 (page 31); power mean is *continuous* and *monontone*: Theorem B.7 (page 31); power triangle function is *continuous* and *monontone*: Corollary 4.6 (page 16); Special values of *p*: Corollary 4.7 (page 17), Corollary B.8 (page 34)

definition using open sets, coincide with each other. Again, this is largely because the open balls of a metric space are open.[5]

However, this is not the case for all *distance space*s. In general, the open balls of a distance space are not open, and they are not a base for a topology. In fact, the open balls of a distance space are a base for a topology if and only if the open balls are open. While the open sets in a distance space do induce a topology, it's open balls may not.[6]

## 2    Standard definitions

### 2.1    Standard sets

**Definition 2.1**    Let $\mathbb{R}$ be the **set of real numbers**. Let $\mathbb{R}^{\vdash} \triangleq \{x \in \mathbb{R} \mid x \geq 0\}$ be the **set of non-negative real numbers**. Let $\mathbb{R}^{+} \triangleq \{x \in \mathbb{R} \mid x > 0\}$ be the **set of postive real numbers**. Let $\mathbb{R}^{*} \triangleq \mathbb{R} \cup \{-\infty, \infty\}$ be the set of **extended real numbers**. [7] Let $\mathbb{Z}$ be the **set of integers**. Let $\mathbb{N} \triangleq \{n \in \mathbb{Z} \mid n \geq 1\}$ be the **set of natural numbers**. Let $\mathbb{Z}^{*} \triangleq \mathbb{Z} \cup \{-\infty, \infty\}$ be the *extended set of integers*.

**Definition 2.2**    Let $X$ be a set. The quantity $2^{X}$ is the *power set of $X$* such that
$$2^{X} \quad \triangleq \quad \{A \subseteq X\} \quad \text{(the set of all subsets of $X$).}$$

### 2.2    Relations

**Definition 2.3**    [8]    Let $X$ and $Y$ be *sets*. The **Cartesian product** $X \times Y$ of $X$ and $Y$ is the set $X \times Y \triangleq \{(x, y) \mid x \in X \text{ and } y \in Y\}$. An **ordered pair** $(x, y)$ on $X$ and $Y$ is any element in $X \times Y$. A **relation** $\circledR$ on $X$ and $Y$ is any subset of $X \times Y$ such that $\circledR \subseteq X \times Y$. The set $2^{XY}$ is the **set of all relations** in $X \times Y$. A *relation* $\mathsf{f} \in 2^{XY}$ is a **function** if $(x, y_1) \in \mathsf{f}$ and $(x, y_2) \in \mathsf{f} \implies y_1 = y_2$. The set $Y^{X}$ is the **set of all functions** in $2^{XY}$.

### 2.3    Set functions

**Definition 2.4**    [9]    Let $2^{X}$ be the *power set* (Definition 2.2 page 4) of a set $X$.
    A set $S(X)$ is a **set structure** on $X$ if        $S(X) \subseteq 2^{X}$.
    A set structure $Q(X)$ is a **paving** on $X$ if        $\varnothing \in Q(X)$.

---

[5] *open ball*: Definition 3.5 page 6; *metric space*: Definition 4.5 page 16; *base*: Definition A.3 page 27; *topology*: Definition A.1 page 26; *open*: Definition 3.6 page 7; *continuity* in *topological space*: Definition A.11 page 28; *convergence* in *distance space*: Definition 3.11 page 9; *convergence* in *topological space*: Definition A.16 page 28; *closed set*: Definition A.1 page 26; *closure, interior, accumulation point*: Definition A.8 page 27; coincidence in all *metric space*s and some *power distance space*s: Theorem 4.15 page 21;

[6] *if and only if* statement: Theorem 3.10 page 8; open sets of a distance space induce a topology: Corollary 3.8 page 8;

[7]  [Rana(2002)] pages 385–388 ⟨Appendix A⟩

[8]  [Maddux(2006)] page 4,  [Halmos(1960)] pages 26–30,  [Suppes(1972)] page 86,  [Kelley(1955)] page 10,  [Bourbaki(1939)],  [Bottazzini(1986)] page 7,  [Comtet(1974)] page 4 ⟨$\mid Y^X \mid$⟩; The notation $Y^X$ and $2^{XY}$ is motivated by the fact that for finite $X$ and $Y$, $\mid Y^X \mid = \mid Y \mid^{\mid X \mid}$ and $\mid 2^{XY} \mid = 2^{\mid X \mid \cdot \mid Y \mid}$.

[9]  [Molchanov(2005)] page 389,  [Pap(1995)] page 7,  [Hahn and Rosenthal(1948)] page 254

**Definition 2.5** [10] Let $\mathcal{Q}(X)$ be a *paving* (Definition 2.4 page 4) on a set $X$. Let $Y$ be a set containing the element 0. A function $\mathsf{m} \in Y^{\mathcal{Q}(X)}$ is a **set function** if $\mathsf{m}(\varnothing) = 0$.

**Definition 2.6** The *set function* (Definition 2.5 page 5) $|A| \in \mathbb{Z}^{*2^X}$ is the *cardinality of A* such that

$$|A| \triangleq \left\{ \begin{array}{ll} \text{the number of elements in } A & \text{for } finite \ A \\ \infty & \text{otherwise} \end{array} \right\} \quad \forall A \in 2^X$$

**Definition 2.7** Let $|X|$ be the *cardinality* (Definition 2.6 page 5) of a *set $X$*.
The structure $\varnothing$ is the **empty set**, and is a *set* such that $|\varnothing| = 0$.

## 2.4 Order

**Definition 2.8** [11] Let $X$ be a set. A relation $\leq$ is an **order relation** in $2^{XX}$ (Definition 2.3 page 4) if

1. $x \leq x$                   $\forall_{x \in X}$   (*reflexive*)    and   **preorder**
2. $x \leq y$  and  $y \leq z \implies x \leq z$  $\forall_{x,y \in X}$  (*transitive*)   and
3. $x \leq y$  and  $y \leq x \implies x = y$  $\forall_{x,y \in X}$  (*anti-symmetric*)

An **ordered set** is the pair $(X, \leq)$.

**Definition 2.9** [12] In an *ordered set* $(X, \leq)$,

the set $[x:y] \triangleq \{z \in X \mid x \leq z \leq y\}$ is a **closed interval**   and
the set $(x:y] \triangleq \{z \in X \mid x < z \leq y\}$ is a **half-open interval**   and
the set $[x:y) \triangleq \{z \in X \mid x \leq z < y\}$ is a **half-open interval**   and
the set $(x:y) \triangleq \{z \in X \mid x < z < y\}$ is an **open interval**.

**Definition 2.10** Let $(\mathbb{R}, \leq)$ be the *ordered set of real numbers* (Definition 2.8 page 5).
The **absolute value** $|\cdot| \in \mathbb{R}^{\mathbb{R}}$ is defined as[13]    $|x| \triangleq \left\{ \begin{array}{ll} -x & \text{for } x \leq 0 \\ x & \text{otherwise} \end{array} \right\}$.

# 3 Background: distance spaces

A *distance space* (Definition 3.1 page 6) can be defined as a *metric space* (Definition 4.5 page 16) without the *triangle inequality* constraint. Much of the material in this section about *distance space*s is standard in *metric space*s. However, this paper works through this material again to demonstrate "how far we can go", and can't go, without the *triangle inequality*.

---

[10] [Hahn and Rosenthal(1948)], 🔖 [Choquet(1954)], ✎ [Pap(1995)] page 8 ⟨Definition 2.3: extended real-valued set function⟩, [Halmos(1950)] page 30 ⟨§7. MEASURE ON RINGS⟩,

[11] [MacLane and Birkhoff(1999)] page 470, [Beran(1985)] page 1, 🔖 [Korselt(1894)] page 156 ⟨I, II, (1)⟩, 🔖 [Dedekind(1900)] page 373 ⟨I–III⟩. An *order relation* is also called a **partial order relation**. An *ordered set* is also called a **partially ordered set** or **poset**.

[12] [Apostol(1975)] page 4, ✎ [Ore(1935)] page 409

[13] A more general definition for *absolute value* is available for any *commutative ring*: Let $R$ be a *commutative ring*. A function $|\cdot|$ in $R^R$ is an **absolute value**, or **modulus**, on $R$ if

1. $|x| \geq 0$                 $x \in \mathbb{R}$   (*non-negative*)    and
2. $|x| = 0 \iff x = 0$  $x \in \mathbb{R}$  (*nondegenerate*)    and
3. $|xy| = |x| \cdot |y|$      $x,y \in \mathbb{R}$  (*homogeneous / submultiplicative*)  and
4. $|x + y| \leq |x| + |y|$   $x,y \in \mathbb{R}$  (*subadditive / triangle inequality*)

Reference: ✎ [Cohn(2002 December 6)] page 312

## 3.1　Fundamental structure of distance spaces

### 3.1.1　Definitions

**Definition 3.1** [14] A function $d$ in the set $\mathbb{R}^{X \times X}$ (Definition 2.3 page 4) is a **distance** if

1. $d(x, y) \geq 0$　　$\forall x, y \in X$　(*non-negative*)　and
2. $d(x, y) = 0 \iff x = y$　$\forall x, y \in X$　(*nondegenerate*)　and
3. $d(x, y) = d(y, x)$　　$\forall x, y \in X$　(*symmetric*)

The pair $(X, d)$ is a **distance space** if $d$ is a *distance* on a set $X$.

**Definition 3.2** [15] Let $(X, d)$ be a *distance space* and $2^X$ be the *power set* of $X$ (Definition 2.2 page 4).　The **diameter** in $(X, d)$ of a set $A \in 2^X$ is　$\operatorname{diam} A \triangleq \begin{cases} 0 & \text{for } A = \varnothing \\ \sup \{d(x, y) \,|\, x, y \in A\} & \text{otherwise} \end{cases}$

**Definition 3.3** [16] Let $(X, d)$ be a *distance space*. Let $2^X$ be the *power set* (Definition 2.2 page 4) of $X$. A set $A$ is **bounded** in $(X, d)$ if $A \in 2^X$ and $\operatorname{diam} A < \infty$.

### 3.1.2　Properties

**Remark 3.4**　Let $(x_n)_{n \in \mathbb{Z}}$ be a *sequence* in a *distance space* $(X, d)$. The *distance space* $(X, d)$ does not necessarily have all the nice properties that a *metric space* (Definition 4.5 page 16) has. In particular, note the following:

1. $d$ is a *distance* in $(X, d)$　　$\implies$　$d$ is *continuous* in $(X, d)$　(Example 3.23 page 14).
2. $B$ is an *open ball* in $(X, d)$　　$\implies$　$B$ is *open* in $(X, d)$　(Example 3.22 page 13).
3. $\boldsymbol{B}$ is the set of all *open ball*s in $(X, d)$　$\implies$　$\boldsymbol{B}$ is a *base* for a topology on $X$　(Example 3.22 page 13).[17]
4. $(x_n)$ is *convergent* in $(X, d)$　　$\implies$　limit is *unique*　(Example 3.21 page 12).
5. $(x_n)$ is *convergent* in $(X, d)$　　$\implies$　$(x_n)$ is *Cauchy* in $(X, d)$　(Example 3.22 page 13).

## 3.2　Open sets in distance spaces

### 3.2.1　Definitions

**Definition 3.5** [18] Let $(X, d)$ be a *distance space* (Definition 3.1 page 6). Let $\mathbb{R}^+$ be the *set of positive real numbers* (Definition 2.1 page 4).

An **open ball** centered at $x$ with radius $r$ is the set　$B(x, r) \triangleq \{y \in X \,|\, d(x, y) < r\}$.

A **closed ball** centered at $x$ with radius $r$ is the set　$\overline{B}(x, r) \triangleq \{y \in X \,|\, d(x, y) \leq r\}$.

---

[14] 📖 [Menger(1928)] page 76 ⟨"Abstand $a$ $b$ definiert ist…" (distance from $a$ to $b$ is defined as…")⟩, 📖 [Wilson(1931)] page 361 ⟨§1., "distance", "semi-metric space"⟩, 📖 [Blumenthal(1938)] page 38, ✎ [Blumenthal(1953)] page 7 ⟨"Definition 5.1. A distance space is called semimetric provided…"⟩, 📖 [Galvin and Shore(1984)] page 67 ⟨"distance function"⟩, ✎ [Laos(1998) page 118] ⟨"distance space"⟩, ✎ [Khamsi and Kirk(2001) page 13] ⟨"semimetric space"⟩, 📖 [Bessenyei and Pales(2014)] page 2 ⟨"semimetric space"⟩, ✎ [Deza and Deza(2014) page 3] ⟨"**distance** (or **dissimilarity**)"⟩

[15] in *metric space*: ✎ [Hausdorff(1937)], page 166, ✎ [Copson(1968)], page 23, ✎ [Michel and Herget(1993)], page 267, ✎ [Molchanov(2005) page 389]

[16] in *metric space*: ✎ [Thron(1966)], page 154 ⟨definition 19.5⟩, ✎ [Bruckner et al.(1997)Bruckner, Bruckner, and Thomson) page 356]

[17] 📖 [Heath(1961)] page 810 ⟨Theorem⟩, 📖 [Galvin and Shore(1984)] page 71 ⟨2.3 Lemma⟩

[18] in *metric space*: ✎ [Aliprantis and Burkinshaw(1998)], page 35

**Definition 3.6** Let $(X, \mathsf{d})$ be a *distance space*. Let $X \setminus A$ be the *set difference* of $X$ and a set $A$. A set $U$ is **open** in $(X, \mathsf{d})$ if $U \in 2^X$ and for every $x$ in $U$ there exists $r \in \mathbb{R}^+$ such that $\mathsf{B}(x, r) \subseteq U$. A set $U$ is an **open set** in $(X, \mathsf{d})$ if $U$ is *open* in $(X, \mathsf{d})$. A set $D$ is **closed** in $(X, \mathsf{d})$ if $(X \setminus D)$ is *open*. A set $D$ is a **closed set** in $(X, \mathsf{d})$ if $D$ is *closed* in $(X, \mathsf{d})$.

### 3.2.2   Properties

**Theorem 3.7** [19] *Let $(X, \mathsf{d})$ be a* distance space. *Let $N$ be any (finite) positive integer. Let $\Gamma$ be a* set *possibly with an uncountable number of elements.*

| | | | |
|---|---|---|---|
| *1.* | | $X$ | *is* OPEN. |
| *2.* | | $\varnothing$ | *is* OPEN. |
| *3.* | *each element in* $\left\{ U_n \big|_{n=1,2,\ldots,N} \right\}$ | *is* OPEN $\implies$ $\displaystyle\bigcap_{n=1}^{N} U_n$ | *is* OPEN. |
| *4.* | *each element in* $\left\{ U_\gamma \in 2^X \mid \gamma \in \Gamma \right\}$ | *is* OPEN $\implies$ $\displaystyle\bigcup_{\gamma \in \Gamma} U_\gamma$ | *is* OPEN. |

✎ PROOF:

(1) Proof that $X$ is *open* in $(X, \mathsf{d})$:

    (a) By definition of *open set* (Definition 3.6 page 7), $X$ is *open* $\iff \forall x \in X \ \ \exists r$ such that $\mathsf{B}(x, r) \subseteq X$.

    (b) By definition of *open ball* (Definition 3.5 page 6), it is always true that $\mathsf{B}(x, r) \subseteq X$ in $(X, \mathsf{d})$.

    (c) Therefore, $X$ is *open* in $(X, \mathsf{d})$.

(2) Proof that $\varnothing$ is *open* in $(X, \mathsf{d})$:

    (a) By definition of *open set* (Definition 3.6 page 7), $\varnothing$ is *open* $\iff \forall x \in X \ \ \exists r$ such that $\mathsf{B}(x, r) \subseteq \varnothing$.

    (b) By definition of *empty set* $\varnothing$ (Definition 2.7 page 5), this is always true because no $x$ is in $\varnothing$.

    (c) Therefore, $\varnothing$ is *open* in $(X, \mathsf{d})$.

(3) Proof that $\bigcup U_\gamma$ is *open* in $(X, \mathsf{d})$:

    (a) By definition of *open set* (Definition 3.6 page 7), $\bigcup U_\gamma$ is *open* $\iff \forall x \in \bigcup U_\gamma \ \ \exists r$ such that $\mathsf{B}(x, r) \subseteq \bigcup U_\gamma$.

    (b) If $x \in \bigcup U_\gamma$, then there is at least one $U \in \bigcup U_\gamma$ that contains $x$.

    (c) By the left hypothesis in (4), that set $U$ is open and so for that $x$, $\exists r$ such that $\mathsf{B}(x, r) \subseteq U \subseteq \bigcup U_\gamma$.

    (d) Therefore, $\bigcup U_\gamma$ is *open* in $(X, \mathsf{d})$.

(4) Proof that $U_1$ and $U_2$ are *open* $\implies U_1 \cap U_2$ is *open*:

    (a) By definition of *open set* (Definition 3.6 page 7), $U_1 \cap U_2$ is *open* $\iff \forall x \in U_1 \cap U_2 \ \ \exists r$ such that $\mathsf{B}(x, r) \subseteq U_1 \cap U_2$.

    (b) By the left hypothesis above, $U_1$ and $U_2$ are *open*; and by the definition of *open set*s (Definition 3.6 page 7), there exists $r_1$ and $r_2$ such that $\mathsf{B}(x, r_1) \subseteq U_1$ and $\mathsf{B}(x, r_2) \subseteq U_2$.

    (c) Let $r \triangleq \min \{r_1, r_2\}$. Then $\mathsf{B}(x, r) \subseteq U_1$ and $\mathsf{B}(x, r) \subseteq U_2$.

    (d) By definition of *set intersection* $\cap$ then, $\mathsf{B}(x, r) \subseteq U_1 \cap U_2$.

    (e) By definition of *open set* (Definition 3.6 page 7), $U_1 \cap U_2$ is *open*.

(5) Proof that $\bigcap_{n=1}^{N} U_n$ is *open* (by induction):

    (a) Proof for $N = 1$ case: $\bigcap_{n=1}^{N} U_n = \bigcap_{n=1}^{1} U_n = U_1$ is *open* by hypothesis.

---

[19] in *metric space*: ✎ [Dieudonné(1969)], pages 33–34, ✎ [Rosenlicht(1968)] page 39

(b)  Proof that $N$ case $\implies$ $N+1$ case:

$$\bigcap_{n=1}^{N+1} U_n = \left(\bigcap_{n=1}^{N} U_n\right) \cap U_{N+1} \qquad \text{by property of } \bigcap$$

$$\implies \quad open \qquad\qquad\qquad \text{by "}N\text{ case" hypothesis and (4) lemma page 7}$$

**Corollary 3.8**  *Let* $(X, \mathsf{d})$ *be a* DISTANCE SPACE. *The set*  $T \triangleq \left\{ U \in 2^X \mid U \text{ is } \text{OPEN } in (X, \mathsf{d}) \right\}$  *is a* TOPOLOGY *on* $X$, *and* $(X, T)$ *is a* TOPOLOGICAL SPACE.

✎PROOF:   This follows directly from the definition of an *open* set (Definition 3.6 page 7), Theorem 3.7 (page 7), and the definition of *topology* (Definition A.1 page 26).

Of course it is possible to define a very large number of topologies even on a finite set with just a handful of elements;[20] and it is possible to define an infinite number of topologies even on a *linearly ordered* infinite set like the *real line* ($\mathbb{R}$, $\leq$).[21] Be that as it may, Definition 3.9 (next definition) defines a single but convenient *topological space* in terms of a *distance space*. Note that every *metric space* conveniently and naturally induces a *topological space* because the *open ball*s of the metric space form a *base* for the *topology*. This is not the case for all distance spaces. But if the open balls of a *distance space* are all *open*, then those open balls induce a topology (next theorem).[22]

**Definition 3.9**  Let $(X, \mathsf{d})$ be a *distance space*. The set  $T \triangleq \left\{ U \in 2^X \mid U \text{ is } open \text{ in } (X, \mathsf{d}) \right\}$  is the **topology induced by** $(X, \mathsf{d})$ **on** $X$. The pair $(X, T)$ is called the **topological space induced by** $(X, \mathsf{d})$.

For any *distance space* $(X, \mathsf{d})$, no matter how strange, there is guaranteed to be at least one *topological space induced by* $(X, \mathsf{d})$—and that is the *indiscrete topological space* (Example A.2 page 26) because for any distance space $(X, \mathsf{d})$, $\varnothing$ and $X$ are *open set*s in $(X, \mathsf{d})$ (Theorem 3.7 page 7).

**Theorem 3.10**  *Let* $B$ *be the set of all* OPEN BALLS *in a* DISTANCE SPACE $(X, \mathsf{d})$.
$$\left\{ every \text{ OPEN BALL } in \text{ } B \text{ is OPEN} \right\} \quad\iff\quad \left\{ B \text{ is a BASE } for \text{ a TOPOLOGY} \right\}$$

---

[20] For a finite set $X$ with $n$ elements, there are 29 topologies on $X$ if $n = 3$, 6942 topologies on $X$ if $n = 5$, and 8,977,053,873,043 (almost 9 trillion) topologies on $X$ if $n = 10$. References: ⌨ [oei(2014)] ⟨http://oeis.org/A000798⟩, ✎ [Brown and Watson(1996)], page 31, ✎ [Comtet(1974)] page 229, ✎ [Comtet(1966)], ✎ [Chatterji(1967)], page 7, ✎ [Evans et al.(1967)Evans, Harary, and Lynn], ✎ [Krishnamurthy(1966)], page 157

[21] For examples of topologies on the real line, see the following: ✎ [Adams and Franzosa(2008)] page 31 ⟨"six topologies on the real line"⟩, ✎ [Salzmann et al.(2007)Salzmann, Grundhöfer, Hähl, and Löwen] pages 64–70 ⟨Weird topologies on the real line⟩, ✎ [Murdeshwar(1990)] page 53 ⟨"often used topologies on the real line"⟩, ✎ [Joshi(1983)] pages 85–91 ⟨§4.2 Examples of Topological Spaces⟩

[22] *metric space*: Definition 4.5 page 16; *open ball*: Definition 3.5 page 6; *base*: Definition A.3 page 27; *topology*: Definition A.1 page 26; not all open balls are open in a distance space: Example 3.21 (page 12) and Example 3.22 (page 13);

✎ PROOF:

| every *open ball* in ***B*** is *open* |
| :-- |

$\implies$ for every $x$ in $B_y \in \boldsymbol{B}$ there exists $r \in \mathbb{R}^+$ such that $\mathrm{B}(x,r) \subseteq B_y$    by definition of *open* (Definition 3.6 page 7)

$\implies$ $\left\{\begin{array}{l}\text{for every } x \in X \text{ and for every } B_y \in \boldsymbol{B} \text{ containing } x, \\ \text{there exists } B_x \in \boldsymbol{B} \text{ such that } \quad x \in B_x \subseteq B_y.\end{array}\right\}$    because $\forall (x,r) \in X \times \mathbb{R}^+,\, \mathrm{B}(x,r) \subseteq X$

$\implies$ ***B*** is a *base* for ***T***    by Theorem A.4 page 27

$\implies$ $\left\{\begin{array}{l}\text{for every } x \in X \text{ and for every } U \subseteq \boldsymbol{T} \text{ containing } x, \\ \text{there exists } B_x \in \boldsymbol{B} \text{ such that } \quad x \in B_x \subseteq U.\end{array}\right\}$    by Theorem A.4 page 27

$\implies$ $\left\{\begin{array}{l}\text{for every } x \in X \text{ and for every } B_y \in \boldsymbol{B} \subseteq \boldsymbol{T} \text{ containing } x, \\ \text{there exists } B_x \in \boldsymbol{B} \text{ such that } \quad x \in B_x \subseteq B_y.\end{array}\right\}$    by definition of *base* (Definition A.3 page 27)

$\implies$ $\left\{\begin{array}{l}\text{for every } x \in B_y \in \boldsymbol{B} \subseteq \boldsymbol{T}, \\ \text{there exists } B_x \in \boldsymbol{B} \text{ such that } \quad x \in B_x \subseteq B_y.\end{array}\right\}$

$\implies$ | every *open ball* in ***B*** is *open* |    by definition of *open* (Definition 3.6 page 7)

## 3.3 Sequences in distance spaces

### 3.3.1 Definitions

**Definition 3.11** [23] Let $(x_n \in X)_{n \in \mathbb{Z}}$ be a *sequence* in a *distance space* $(X, \mathsf{d})$. The sequence $(x_n)$ **converges** to a **limit** $x$ if for any $\varepsilon \in \mathbb{R}^+$, there exists $N \in \mathbb{Z}$ such that for all $n > N$, $\mathsf{d}(x_n, x) < \varepsilon$.
This condition can be expressed in any of the following forms:

1. The **limit** of the sequence $(x_n)$ is $x$.
2. The sequence $(x_n)$ is **convergent** with limit $x$.
3. $\lim\limits_{n \to \infty} (x_n) = x$.
4. $(x_n) \to x$.

A *sequence* that converges is **convergent**.

**Definition 3.12** [24] Let $(x_n \in X)_{n \in \mathbb{Z}}$ be a *sequence* in a *distance space* $(X, \mathsf{d})$.
The sequence $(x_n)$ is a **Cauchy sequence** in $(X, \mathsf{d})$ if

     for every $\varepsilon \in \mathbb{R}^+$, there exists $N \in \mathbb{Z}$ such that $\forall n, m > N$, $\mathsf{d}(x_n, x_m) < \varepsilon$    (*Cauchy condition*).

**Definition 3.13** [25] Let $(x_n \in X)_{n \in \mathbb{Z}}$ be a *sequence* in a *distance space* $(X, \mathsf{d})$.
The sequence $(x_n \in X)_{n \in \mathbb{Z}}$ is **complete** in $(X, \mathsf{d})$ if

     $(x_n)$ is *Cauchy* in $(X, \mathsf{d})$    $\implies$    $(x_n)$ is *convergent* in $(X, \mathsf{d})$.

### 3.3.2 Properties

**Proposition 3.14** [26] *Let* $(x_n \in X)_{n \in \mathbb{Z}}$ *be a* SEQUENCE *in a* DISTANCE SPACE $(X, \mathsf{d})$.

     $\big\{ (x_n) \text{ is CAUCHY } in (X, \mathsf{d}) \big\}$    $\implies$    $\big\{ (x_n) \text{ is BOUNDED } in (X, \mathsf{d}) \big\}$

✎Proof:

$(x_n)$ is *Cauchy* $\implies$ for every $\varepsilon \in \mathbb{R}^+$, $\exists N \in \mathbb{Z}$  such that  $\forall n, m > N$, $\mathrm{d}(x_n, x_m) < \varepsilon$     (by Definition 3.12 page 9)

$\implies \exists N \in \mathbb{Z}$  such that  $\forall n, m > N$, $\mathrm{d}(x_n, x_m) < 1$     (arbitrarily choose $\varepsilon \triangleq 1$)

$\implies \exists N \in \mathbb{Z}$  such that  $\forall n, m \in \mathbb{Z}$, $\mathrm{d}(x_n, x_{m+1}) < \max\left\{ \{1\} \cup \{\mathrm{d}(x_p, x_q) \,|\, p, q \not> N\}\right\}$

$\implies (x_n)$ is *bounded*     (by Definition 3.3 page 6)

**Proposition 3.15** [27] *Let* $(x_n \in X)_{n \in \mathbb{Z}}$ *be a* SEQUENCE *in a* DISTANCE SPACE $(X, \mathrm{d})$. *Let* $\mathrm{f} \in \mathbb{Z}^{\mathbb{Z}}$ *(Definition 2.3 page 4) be a* STRICTLY MONOTONE *function such that* $\mathrm{f}(n) < \mathrm{f}(n+1)$.

$$\underbrace{(x_n)_{n \in \mathbb{Z}} \text{ is CAUCHY}}_{\text{sequence is CAUCHY}} \implies \underbrace{(x_{\mathrm{f}(n)})_{n \in \mathbb{Z}} \text{ is CAUCHY}}_{\text{subsequence is also CAUCHY}}$$

✎Proof:

$(x_n)_{n \in \mathbb{Z}}$ is *Cauchy*

$\implies$ for any given $\varepsilon > 0$, $\exists N$  such that  $\forall n, m > N$, $\mathrm{d}(x_n, x_m) < \varepsilon$     by Definition 3.12 page 9

$\implies$ for any given $\varepsilon > 0$, $\exists N'$  such that  $\forall \mathrm{f}(n), \mathrm{f}(m) > N'$, $\mathrm{d}(x_{\mathrm{f}(n)}, x_{\mathrm{f}(m)}) < \varepsilon$

$\implies (x_{\mathrm{f}(n)})_{n \in \mathbb{Z}}$ is *Cauchy*     by Definition 3.12 page 9

**Theorem 3.16** [28] *Let* $(X, \mathrm{d})$ *be a* DISTANCE SPACE. *Let* $A^-$ *be the* CLOSURE *(Definition A.8 page 27) of a* $A$ *in a* TOPOLOGICAL SPACE INDUCED *by* $(X, \mathrm{d})$.

$$\left.\begin{array}{ll} 1. & \text{LIMIT} s \text{ } are \text{ UNIQUE } in (X, \mathrm{d}) \quad \text{\scriptsize(Definition 3.11 page 9)} \quad and \\ 2. & (A, \mathrm{d}) \text{ } is \text{ COMPLETE } in (X, \mathrm{d}) \quad \text{\scriptsize(Definition 3.13 page 9)} \end{array}\right\} \implies \underbrace{A \text{ } is \text{ CLOSED } in (X, \mathrm{d})}_{A = A^-}$$

✎Proof:

(1) Proof that $A \subseteq A^-$: by Lemma A.10 page 27

(2) Proof that $A^- \subseteq A$ (proof that $x \in A^- \implies x \in A$):

    (a) Let $x$ be a point in $A^-$ ($x \in A^-$).

    (b) Define a *sequence* of open balls $\left( \mathrm{B}(x, \frac{1}{1}), \mathrm{B}(x, \frac{1}{2}), \mathrm{B}(x, \frac{1}{3}), \ldots \right)$.

    (c) Define a *sequence* of points $(x_1, x_2, x_3, \ldots)$ such that $x_n \in \mathrm{B}(x_n, \frac{1}{n}) \cap A$.

    (d) Then $(x_n)$ is *convergent* in $X$ with limit $x$ by Definition 3.11 page 9

    (e) and $(x_n)$ is *Cauchy* in $A$ by Definition 3.12 page 9.

    (f) By the hypothesis 2, $(x_n)$ is therefore also *convergent* in $A$.
    Let this limit be $y$. Note that $y \in A$.

    (g) By hypothesis 1, limits are *unique*, so $y = x$.

    (h) Because $y \in A$ (item (2f)) and $y = x$ (item (2g)), so $x \in A$.

    (i) Therefore, $x \in A^- \implies x \in A$ and $A^- \subseteq A$.

**Proposition 3.17** [29] *Let* $(x_n)_{n \in \mathbb{Z}}$ *be a sequence in a* DISTANCE SPACE $(X, \mathrm{d})$. *Let* $\mathrm{f} : \mathbb{Z} \to \mathbb{Z}$ *be a strictly*

---

[27] in *metric space*: ✎ [Rosenlicht(1968)] page 52
[28] in *metric space*: ✎ [Kubrusly(2001)] page 128 ⟨Theorem 3.40⟩, ✎ [Haaser and Sullivan(1991)] page 75 ⟨6·10, 6·11 Propositions⟩, ✎ [Bryant(1985)] page 40 ⟨Theorem 3.6, 3.7⟩, ✎ [Sutherland(1975)] pages 123–124
[29] in *metric space*: ✎ [Rosenlicht(1968)] page 46

*increasing function such that* $f(n) < f(n+1)$.

$$\underbrace{(x_n)_{n \in \mathbb{Z}} \to x}_{\textit{sequence converges to limit } x} \quad \implies \quad \underbrace{\left(x_{f(n)}\right)_{n \in \mathbb{Z}} \to x}_{\textit{subsequence converges to the same limit } x}$$

✎Proof:

$$
\begin{aligned}
(x_n)_{n \in \mathbb{Z}} \to x &\implies \forall \varepsilon > 0, \; \exists N \quad \text{such that} \quad \forall n > N, \; d(x_n, x) < \varepsilon && \text{by Theorem 4.15 page 21} \\
&\implies \forall \varepsilon > 0, \; \exists f(N) \quad \text{such that} \quad \forall f(n) > f(N), \; d\left(x_{f(n)}, x\right) < \varepsilon \\
&\implies \left(x_{f(n)}\right)_{n \in \mathbb{Z}} \to x && \text{by Theorem 4.15 page 21}
\end{aligned}
$$

---

**Theorem 3.18**   (Cantor intersection theorem) [30]   *Let* $(X, d)$ *be a* DISTANCE SPACE *(Definition 3.1 page 6)*, $(A_n)_{n \in \mathbb{Z}}$ *a* SEQUENCE *with each* $A_n \in 2^X$, *and* $|A|$ *the number of elements in* $A$.

$$
\left\{
\begin{array}{lll}
1. \; (X, d) \textit{ is } \text{COMPLETE} & \textit{(Definition 3.13 page 9)} & \textit{and} \\
2. \; A_n \textit{ is } \text{CLOSED} & \forall n \in \mathbb{N} \quad \textit{(Definition A.1 page 26)} & \textit{and} \\
3. \; \text{diam } A_n \geq \text{diam } A_{n+1} & \forall n \in \mathbb{N} \quad \textit{(Definition 3.2 page 6)} & \textit{and} \\
4. \; \text{diam } (A_n)_{n \in \mathbb{Z}} \to 0 & \textit{(Definition 3.11 page 9)}
\end{array}
\right\}
\implies
\left\{ \left| \bigcap_{n \in \mathbb{N}} A_n \right| = 1 \right\}
$$

✎Proof:

(1) Proof that $\left| \bigcap_{n \in \mathbb{Z}} A_n \right| < 2$:

    (a) Let $A \triangleq \cap A_n$.

    (b) $x \neq y$ and $\{x, y\} \in A \implies d(x, y) > 0$ and $\{x, y\} \subseteq A_n \forall n$

    (c) $\exists n$ such that $\text{diam } A_n < d(x, y)$ by left hypothesis 4

    (d) $\implies \exists n$ such that $\sup \left\{ d(x, y) \, \middle| \, x, y \in A_n \right\} < d(x, y)$

    (e) This is a contradiction, so $\{x, y\} \notin A$ and $\left| \bigcap A_n \right| < 2$.

(2) Proof that $\left| \cap A_n \right| \geq 1$:

    (a) Let $x_n \in A_n$ and $x_m \in A_m$

    (b) $\forall \varepsilon, \; \exists N \in \mathbb{N}$ such that $A_N < \varepsilon$

    (c) $\forall m, n > N, \; x_n \in A_n \subseteq A_N$ and $x_m \in A_m \subseteq A_N$

    (d) $d(x_n, x_m) \leq \text{diam } A_N < \varepsilon \implies \{x_n\}$ is a Cauchy sequence

    (e) Because $\{x_n\}$ is complete, $x_n \to x$.

    (f) $\implies x \in (A_n)^- = A_n$

    (g) $\implies |A_n| \geq 1$

---

**Definition 3.19**   [31] Let $(X, d)$ be a *distance space*. Let $C$ be the set of all *convergent* sequences in $(X, d)$. The *distance function* d is **continuous** in $(X, d)$ if

$$(x_n), (y_n) \in C \implies \lim_{n \to \infty} \left( d(x_n, y_n) \right) = d\left( \lim_{n \to \infty} (x_n), \lim_{n \to \infty} (y_n) \right).$$

A *distance function* is **discontinuous** if it is not *continuous*.

---

[30] in *metric space*: ✎ [Davis(2005)], page 28, ✎ [Hausdorff(1937)], page 150
[31] ✎ [Blumenthal(1953)] page 9 ⟨DEFINITION 6.3⟩

**Remark 3.20**  Rather than defining *continuity* of a *distance function* in terms of the *sequential characterization of continuity* as in Definition 3.19 (previous), we could define continuity using an *inverse image characterization of continuity*" (Definition 3.9 page 8). Assuming an equivalent *topological space* is used for both characterizations, the two characterizations are equivalent (Theorem A.20 page 29). In fact, one could construct an equivalence such as the following:

$$\left\{\begin{array}{l} \mathsf{d} \text{ is } continuous \text{ in } \mathbb{R}^{X^2} \\ \text{(Definition A.11 page 28)} \\ \text{(inverse image characterization of continuity)} \end{array}\right\} \iff \left\{\begin{array}{l} \left(\!\left(x_n\right)\!\right),\left(\!\left(y_n\right)\!\right) \in C \implies \\ \lim_{n\to\infty}\left(\mathsf{d}\left(x_n,y_n\right)\right) = \mathsf{d}\left(\lim_{n\to\infty}\left(\!\left(x_n\right)\!\right),\lim_{n\to\infty}\left(\!\left(y_n\right)\!\right)\right) \\ \text{(Definition A.16 page 28)} \\ \text{(sequential characterization of continuity)} \end{array}\right\}$$

Note that just as $\left(\!\left(x_n\right)\!\right)$ is a sequence in $X$, so the ordered pair $\left(\left(\!\left(x_n\right)\!\right),\left(\!\left(y_n\right)\!\right)\right)$ is a sequence in $X^2$. The remainder follows from Theorem A.20 (page 29). However, use of the *inverse image characterization* is somewhat troublesome because we would need a topology on $X^2$, and we don't immediately have one defined and ready to use. In fact, we don't even immediately have a distance space on $X^2$ defined or even open balls in such a distance space. The result is, for the scope of this paper, it is arguably not worthwhile constructing the extra structure, but rather instead this paper uses the *sequential characterization* as a definition (as in Definition 3.19).

## 3.4   Examples

Similar distance functions and several of the observations for the examples in this section can be found in ✎ [Blumenthal(1953)] pages 8–13.

In a *metric space*, all *open ball*s are *open*, the *open ball*s form a *base* for a *topology*, the limits of *convergent* sequences are *unique*, and the *metric function* is *continuous*. In the *distance space* of the next example, none of these properties hold.

**Example 3.21**  [32] Let $(x,y)$ be an *ordered pair* in $\mathbb{R}^2$. Let $(a:b)$ be an *open interval* and $(a:b]$ a *half-open interval* in $\mathbb{R}$. Let $|x|$ be the *absolute value* of $x\in\mathbb{R}$. The function $\mathsf{d}(x,y)\in\mathbb{R}^{\mathbb{R}\times\mathbb{R}}$ such that

$$\mathsf{d}(x,y) \triangleq \left\{\begin{array}{ll} y & \forall(x,y)\in\{4\}\times(0:2] & \text{(vertical half-open interval)} \\ x & \forall(x,y)\in(0:2]\times\{4\} & \text{(horizontal half-open interval)} \\ |x-y| & \text{otherwise} & \text{(Euclidean)} \end{array}\right\} \text{ is a } distance \text{ on } \mathbb{R}.$$

Note some characteristics of the *distance space* $(\mathbb{R},\mathsf{d})$:

(1)  $(\mathbb{R},\mathsf{d})$ is not a *metric space* because $\mathsf{d}$ does not satisfy the *triangle inequality*:
$$\mathsf{d}(0,4) \triangleq |0-4| = 4 \nleq 2 = |0-1| + 1 \triangleq \mathsf{d}(0,1) + \mathsf{d}(1,4)$$

(2)  Not every *open ball* in $(\mathbb{R},\mathsf{d})$ is *open*.
For example, the *open ball* $\mathsf{B}(3,2)$ is *not open* because $4\in\mathsf{B}(3,2)$ *but* for all $0<\varepsilon<1$
$$\mathsf{B}(4,\varepsilon) = (4-\varepsilon:4+\varepsilon)\cup(0:\varepsilon) \nsubseteq (1:5) = \mathsf{B}(3,2)$$

(3)  The *open ball*s of $(\mathbb{R},\mathsf{d})$ do not form a *base* for a *topology* on $\mathbb{R}$.
This follows directly from item (2) and Theorem 3.10 (page 8).

[32] A similar distance function $\mathsf{d}$ and item (4) page 13 can in essence be found in ✎ [Blumenthal(1953)] page 8. Definitions for Example 3.21: $(x,y)$: Definition 2.3 (page 4); $(a:b)$ and $(a:b]$: Definition 2.9 (page 5); $|x|$: Definition 2.10 (page 5); $\mathbb{R}^{\mathbb{R}\times\mathbb{R}}$: Definition 2.3 (page 4); *distance*: Definition 3.1 (page 6); *open ball*: Definition 3.5 (page 6); *open*: Definition 3.6 (page 7); *base*: Definition A.3 (page 27); *topology*: Definition A.1 (page 26); *open set*: Definition 3.6 (page 7); *topological space induced by* $(\mathbb{R},\mathsf{d})$: Definition 3.9 (page 8); *discontinuous*: Definition 3.19 (page 11);

(4)  In the *distance space* $(\mathbb{R}, \mathsf{d})$, limits are *not unique*;
For example, the sequence $(1/n)_1^\infty$ converges both to the limit 0 and the limit 4 in $(\mathbb{R}, \mathsf{d})$:

$$\lim_{n\to\infty} \mathsf{d}(x_n, 0) \triangleq \lim_{n\to\infty} \mathsf{d}(1/n, 0) \triangleq \lim_{n\to\infty} |1/n - 0| = 0 \implies (1/n) \to 0$$

$$\lim_{n\to\infty} \mathsf{d}(x_n, 4) \triangleq \lim_{n\to\infty} \mathsf{d}(1/n, 4) \triangleq \lim_{n\to\infty} (1/n) = 0 \implies (1/n) \to 4$$

(5)  The *topological space* $(X, \boldsymbol{T})$ *induced by* $(\mathbb{R}, \mathsf{d})$ also yields limits of 0 and 4 for the sequence $(1/n)_1^\infty$, just as it does in item (4). This is largely due to the fact that, for small $\varepsilon$, the open balls $\mathsf{B}(0, \varepsilon)$ and $\mathsf{B}(4, \varepsilon)$ are *open*.

$\mathsf{B}(0, \varepsilon)$ is *open* $\implies$ for each $U \in \boldsymbol{T}$ that contains 0, $\exists N \in \mathbb{N}$   such that   $1/n \in U$   $\forall n > N$

$\iff (1/n) \to 0$   by definition of *convergence* (Definition A.16 page 28)

$\mathsf{B}(4, \varepsilon)$ is *open* $\implies$ for each $U \in \boldsymbol{T}$ that contains 4, $\exists N \in \mathbb{N}$   such that   $1/n \in U$   $\forall n > N$

$\iff (1/n) \to 4$   by definition of *convergence* (Definition A.16 page 28)

(6)  The distance function $\mathsf{d}$ is *discontinuous* (Definition 3.19 page 11):

$$\lim_{n\to\infty} (\mathsf{d}(1 - 1/n, 4 - 1/n)) = \lim_{n\to\infty} (|(1 - 1/n) - (4 - 1/n)|) = |1 - 4| = 3 \neq 4 = \mathsf{d}(0, 4)$$

$$= \mathsf{d}\left(\lim_{n\to\infty} (1 - 1/n), \lim_{n\to\infty} (4 - 1/n)\right)$$

In a *metric space*, all *convergent* sequences are also *Cauchy*. However, this is not the case for all *distance space*s, as demonstrated next:

**Example 3.22**  [33] The function $\mathsf{d}(x, y) \in \mathbb{R}^{\mathbb{R} \times \mathbb{R}}$ such that

$$\mathsf{d}(x, y) \triangleq \left\{ \begin{array}{ll} |x - y| & \text{for } x = 0 \text{ or } y = 0 \text{ or } x = y \quad \textit{(Euclidean)} \\ 1 & \text{otherwise} \hspace{3.2cm} \textit{(discrete)} \end{array} \right\} \quad \text{is a } \textit{distance} \text{ on } \mathbb{R}.$$

Note some characteristics of the *distance space* $(\mathbb{R}, \mathsf{d})$:

(1)  $(X, \mathsf{d})$ is not a *metric space* because the *triangle inequality* does not hold:
$$\mathsf{d}\left(\tfrac{1}{4}, \tfrac{1}{2}\right) = 1 \not\leq \tfrac{3}{4} = \left|\tfrac{1}{4} - 0\right| + \left|0 - \tfrac{1}{2}\right| = \mathsf{d}\left(\tfrac{1}{4}, 0\right) + \mathsf{d}\left(0, \tfrac{1}{2}\right)$$

(2)  The *open ball* $\mathsf{B}\left(\tfrac{1}{4}, \tfrac{1}{2}\right)$ is *not open* because for any $\varepsilon \in \mathbb{R}^+$, no matter how small,
$$\mathsf{B}(0, \varepsilon) = (-\varepsilon \,:\, +\varepsilon) \not\subseteq \left\{0, \tfrac{1}{3}\right\} = \left\{x \in X \,\middle|\, \mathsf{d}\left(\tfrac{1}{4}, x\right) < \tfrac{1}{2}\right\} \triangleq \mathsf{B}\left(\tfrac{1}{4}, \tfrac{1}{2}\right)$$

(3)  Even though not all the *open ball*s are *open*, it is still possible to have an *open set* in $(X, \mathsf{d})$. For example, the set $U \triangleq \{1, 2\}$ is *open*:
$$\mathsf{B}(1, 1) \triangleq \{x \in X \,|\, \mathsf{d}(1, x) < 1\} = \{1\} \subseteq \{1, 2\} \triangleq U$$
$$\mathsf{B}(2, 1) \triangleq \{x \in X \,|\, \mathsf{d}(2, x) < 1\} = \{2\} \subseteq \{1, 2\} \triangleq U$$

(4)  By item (2) and Theorem 3.10 (page 8), the *open ball*s of $(\mathbb{R}, \mathsf{d})$ do not form a *base* for a *topology* on $\mathbb{R}$.

(5)  Even though the open balls in $(\mathbb{R}, \mathsf{d})$ do not induce a topology on $X$, it is still possible to find a set of *open set*s in $(X, \mathsf{d})$ that *is* a topology. For example, the set $\{\emptyset, \{1, 2\}, \mathbb{R}\}$ is a topology on $\mathbb{R}$.

(6)  In $(\mathbb{R}, \mathsf{d})$, limits of *convergent* sequences are *unique*:
$$(x_n) \to x \implies \lim_{n\to\infty} \mathsf{d}(x_n, x) = \left\{ \begin{array}{lll} \lim |x_n - 0| & = 0 & \text{for } x = 0 \hspace{1.8cm} \text{OR} \\ |x - x| & = 0 & \text{for constant } (x_n) \text{ for } n > N \quad \text{OR} \\ 1 & \neq 0 & \text{otherwise} \end{array} \right\}$$

[33] The distance function $\mathsf{d}$ and item (7) page 14 can in essence be found in 📖 [Blumenthal(1953)] page 9

which says that there are only two ways for a sequence to converge: either $x = 0$ or the sequence eventually becomes constant (or both). Any other sequence will *diverge*. Therefore we can say the following:

  (a) If $x = 0$ and the sequence is not constant, then the limit is *unique* and $0$.

  (b) If $x = 0$ and the sequence is constant, then the limit is *unique* and $0$.

  (c) If $x \neq 0$ and the sequence is constant, then the limit is *unique* and $x$.

  (d) If $x \neq 0$ and the sequence is not constant, then the sequence diverges and there is no limit.

(7) In $(\mathbb{R}, \mathsf{d})$, a *convergent* sequence is not necessarily *Cauchy*. For example,

  (a) the sequence $(1/n)_{n \in \mathbb{N}}$ is *convergent* with limit 0: $\lim_{n \to \infty} \mathsf{d}(1/n, 0) = \lim_{n \to \infty} 1/n = 0$

  (b) However, even though $(1/n)$ is *convergent*, it is *not Cauchy*: $\lim_{n,m \to \infty} \mathsf{d}(1/n, 1/m) = 1 \neq 0$

(8) The *distance function* $\mathsf{d}$ is *discontinuous* in $(X, \mathsf{d})$:

$$\lim_{n \to \infty} \left( \mathsf{d}(1/n, 2 - 1/n) \right) = 1 \neq 2 = \mathsf{d}(0, 2) = \mathsf{d}\left( \lim_{n \to \infty} (1/n), \lim_{n \to \infty} (2 - 1/n) \right).$$

**Example 3.23** [34] The function $\mathsf{d}(x, y) \in \mathbb{R}^{\mathbb{R} \times \mathbb{R}}$ such that

$$\mathsf{d}(x, y) \triangleq \left\{ \begin{array}{ll} 2|x - y| & \forall (x, y) \in \{(0, 1), (1, 0)\} \quad \text{\textit{(dilated Euclidean)}} \\ |x - y| & \text{otherwise} \quad\quad\quad\quad\quad \text{\textit{(Euclidean)}} \end{array} \right\} \quad \text{is a \textit{distance} on } \mathbb{R}.$$

Note some characteristics of the *distance space* $(\mathbb{R}, \mathsf{d})$:

(1) $(\mathbb{R}, \mathsf{d})$ is *not a metric space* because $\mathsf{d}$ does *not* satisfy the *triangle inequality*:
$\mathsf{d}(0, 1) \triangleq 2|0 - 1| = 2 \not\leq 1 = |0 - 1/2| + |1/2 - 1| \triangleq \mathsf{d}(0, 1/2) + \mathsf{d}(1/2, 1)$

(2) The function $\mathsf{d}$ is *discontinuous*:
$$\lim_{n \to \infty} \left( \mathsf{d}(1 - 1/n, 1/n) \right) = \lim_{n \to \infty} \left( |1 - 1/n - 1/n| \right) = 1 \neq 2 = 2|0 - 1| \triangleq \mathsf{d}(0, 1) = \mathsf{d}\left( \lim_{n \to \infty} (1 - 1/n), \lim_{n \to \infty} (1/n) \right).$$

(3) In $(X, \mathsf{d})$, *open ball*s are *open*:

  (a) $\mathsf{p}(x, y) \triangleq |x - y|$ is a *metric* and thus all open balls in that do not contain both 0 and 1 are *open*.

  (b) By Example C.4 (page 36), $\mathsf{q}(x, y) \triangleq 2|x - y|$ is also a *metric* and thus all open balls containing 0 and 1 only are *open*.

  (c) The only question remaining is with regards to open balls that contain 0, 1 and some other element(s) in $\mathbb{R}$. But even in this case, open balls are still open. For example:
$\mathsf{B}(-1, 2) = (-1 : 2) = (-1 : 1) \cup (1 : 2)$
Note that both $(-1 : 1)$ and $(1 : 2)$ are *open*, and thus by Theorem 3.7 (page 7), $\mathsf{B}(-1, 2)$ is *open* as well.

(4) By item (3) and Theorem 3.10 (page 8), the *open ball*s of $(\mathbb{R}, \mathsf{d})$ *do* form a *base* for a *topology* on $\mathbb{R}$.

(5) In $(X, \mathsf{d})$, the limits of *convergent* sequences are *unique*. This is demonstrated in Example 4.22 (page 25) using additional structure developed in Section 4.

(6) In $(X, \mathsf{d})$, *convergent* sequences are *Cauchy*. This is also demonstrated in Example 4.22 (page 25).

---

[34] The distance function $\mathsf{d}$ and item (2) page 14 can in essence be found in  [Blumenthal(1953)] page 9

The *distance functions* in Example 3.21 (page 12)–Example 3.23 (page 14) were all *discontinuous*. In the absence of the *triangle inequality* and in light of these examples, one might try replacing the *triangle inequality* with the weaker requirement of *continuity*. However, as demonstrated by the next example, this also leads to an arguably disastrous result.

**Example 3.24** [35] The function $d \in \mathbb{R}^{\mathbb{R} \times \mathbb{R}}$ such that $d(x, y) \triangleq (x - y)^2$ is a *distance* on $\mathbb{R}$. Note some characteristics of the *distance space* $(\mathbb{R}, d)$:

(1) $(\mathbb{R}, d)$ is *not* a *metric space* because the *triangle inequality* does not hold:
$$d(0, 2) \triangleq (0 - 2)^2 = 4 \nleq 2 = (0 - 1)^2 + (1 - 2)^2 \triangleq d(0, 1) + d(1, 2)$$

(2) The *distance function* $d$ is *continuous* in $(X, d)$. This is demonstrated in the more general setting of Section 4 in Example 4.23 (page 25).

(3) Calculating the length of curves in $(X, d)$ leads to a paradox:[36]

  (a) Partition $[0 : 1]$ into $2^N$ consecutive line segments connected at the points
  $$\left(0, \frac{1}{2^N}, \frac{2}{2^N}, \frac{3}{2^N}, \ldots, \frac{2^{N-1}}{2^N}, 1\right)$$

  (b) Then the distance, as measured by $d$, between any two consecutive points is
  $$d(p_n, p_{n+1}) \triangleq (p_n - p_{n+1})^2 = \left(\frac{1}{2^N}\right)^2 = \frac{1}{2^{2N}}$$

  (c) But this leads to the paradox that the total length of $[0 : 1]$ is 0:
  $$\lim_{N \to \infty} \sum_{n=0}^{2^N - 1} \frac{1}{2^{2N}} = \lim_{N \to \infty} \frac{2^N}{2^{2N}} = \lim_{N \to \infty} \frac{1}{2^N} = 0$$

# 4 Distance spaces with power triangle inequalities

## 4.1 Definitions

This paper introduces a new relation called the *power triangle inequality* (Definition 4.3 page 16). It is a generalization of other common relations, including the *triangle inequality* (Definition 4.4 page 16). The *power triangle inequality* is defined in terms of a function herein called the *power triangle function* (next definition). This function is a special case of the *power mean* with $N = 2$ and $\lambda_1 = \lambda_2 = \frac{1}{2}$ (Definition B.6 page 31). *Power mean*s have the attractive properties of being *continuous* and *strictly monontone* with respect to a free parameter $p \in \mathbb{R}^*$ (Theorem B.7 page 31). This fact is inherited and exploited by the *power triangle inequality* (Corollary 4.6 page 16).

**Definition 4.1** Let $(X, d)$ be a *distance space* (Definition 3.1 page 6). Let $\mathbb{R}^+$ be the set of all *positive real numbers* and $\mathbb{R}^*$ be the set of *extended real numbers* (Definition 2.1 page 4). The **power triangle function** $\tau$ on $(X, d)$ is defined as

$$\tau(p, \sigma; x, y, z; d) \triangleq 2\sigma \left[\frac{1}{2} d^p(x, z) + \frac{1}{2} d^p(z, y)\right]^{\frac{1}{p}} \quad \forall_{(p,\sigma) \in \mathbb{R}^* \times \mathbb{R}^+}, \quad x, y, z \in X$$

---

[35]  [Blumenthal(1953)] pages 12–13,  [Laos(1998)] pages 118–119

[36] This is the method of "inscribed polygons" for calculating the length of a curve and goes back to Archimedes:  [Brunschwig et al.(2003)Brunschwig, Lloyd, and Pellegrin] page 26,  [Walmsley(1920)], page 200 ⟨§158⟩,

**Remark 4.2** [37] In the field of *probabilistic metric spaces*, a function called he *triangle function* was introduced by Sherstnev in 1962. However, the *power triangle function* as defined in this present paper is *not* a special case of (is not compatible with) the *triangle function* of Sherstnev. Another definition of *triangle function* has been offered by Bessenyei in 2014 with special cases of $\Phi(u,v) \triangleq c(u+v)$ and $\Phi(u,v) \triangleq (u^p + v^p)^{\frac{1}{p}}$, which *are* similar to the definition of *power triangle function* offered in this present paper.

**Definition 4.3**   Let $(X,\mathsf{d})$ be a *distance space*. Let $2^{XXX}$ be the set of all trinomial *relation*s (Definition 2.3 page 4) on $X$. A relation $⦻(p,\sigma;\mathsf{d})$ in $2^{XXX}$ is a **power triangle inequality** on $(X,\mathsf{d})$ if
$$⦻(p,\sigma;\mathsf{d}) \triangleq \big\{ (x,y,z) \in X^3 \,|\, \mathsf{d}(x,y) \le \tau(p,\sigma;x,y,z;\mathsf{d}) \big\} \qquad \text{for some } (p,\sigma) \in \mathbb{R}^* \times \mathbb{R}^+.$$
The tupple $(X,\mathsf{d},p,\sigma)$ is a **power distance space** and $\mathsf{d}$ a **power distance** or **power distance function** if $(X,\mathsf{d})$ is a *distance space* in which the *triangle relation* $⦻(p,\sigma;\mathsf{d})$ holds.

The *power triangle function* can be used to define some standard inequalities (next definition). See Corollary 4.7 (page 17) for some justification of the definitions.

**Definition 4.4**   [38] Let $⦻(p,\sigma;\mathsf{d})$ be a *power triangle inequality* on a *distance space* $(X,\mathsf{d})$.

1. $⦻(\infty, ^1\!/_2; \mathsf{d}\,)$ is the **σ-inframetric inequality**
2. $⦻(\infty, \frac{1}{2}; \mathsf{d}\,)$ is the **inframetric inequality**
3. $⦻(2, \sqrt{2}/_2; \mathsf{d}\,)$ is the **quadratic inequality**
4. $⦻(1, \sigma; \mathsf{d}\,)$ is the **relaxed triangle inequality**
5. $⦻(1, 1; \mathsf{d}\,)$ is the **triangle inequality**
6. $⦻(-^1\!/_2, 2; \mathsf{d}\,)$ is the **square mean root inequality**
7. $⦻(0, \frac{1}{3}; \mathsf{d}\,)$ is the **geometric inequality**
8. $⦻(-1, \frac{1}{4}; \mathsf{d}\,)$ is the **harmonic inequality**
9. $⦻(-\infty, \frac{1}{2}; \mathsf{d}\,)$ is the **minimal inequality**

**Definition 4.5**   [39] Let $(X,\mathsf{d})$ be a *distance space* (Definition 3.1 page 6).

1. $(X,\mathsf{d})$ is a      **metric space**          if the *triangle inequality*          holds in $X$.
2. $(X,\mathsf{d})$ is a      **near metric space**      if the *relaxed triangle inequality*   holds in $X$.
3. $(X,\mathsf{d})$ is an     **inframetric space**      if the *inframetric inequality*        holds in $X$.
4. $(X,\mathsf{d})$ is a      **σ-inframetric space**    if the *σ-inframetric inequality*      holds in $X$.

## 4.2   Properties

### 4.2.1   Relationships of the power triangle function

**Corollary 4.6**   *Let* $\tau(p,\sigma;x,y,z;\mathsf{d})$ *be the* POWER TRIANGLE FUNCTION *(Definition 4.1 page 15) in the* DISTANCE SPACE *(Definition 3.1 page 6)* $(X,\mathsf{d})$. *Let* $(\mathbb{R},|\cdot|,\le)$ *be the* ORDERED METRIC SPACE *with the usual ordering relation*

---

[37]  [Sherstnev(1962)], page 4, ✎ [Schweizer and Sklar(1983)] page 9 ⟨(1.6.1)–(1.6.4)⟩, ▣ [Bessenyei and Pales(2014)] page 2

[38]  [Bessenyei and Pales(2014)] page 2, ▣ [Czerwik(1993)] page 5 ⟨*b-metric*; (1),(2),(5)⟩, ▣ [Fagin et al.(2003a)Fagin, Kumar, and Sivakumar], ▣ [Fagin et al.(2003b)Fagin, Kumar, and Sivakumar] ⟨Definition 4.2 (Relaxed metrics)⟩, ▣ [Xia(2009)] page 453 ⟨Definition 2.1⟩, ▣ [Heinonen(2001)] page 109 ⟨14.1 Quasi-metric spaces.⟩, ✎ [Kirk and Shahzad(2014)] page 113 ⟨Definition 12.1⟩, ✎ [Deza and Deza(2014)] page 7, ▣ [Hoehn and Niven(1985)] page 151, ✎ [Gibbons et al.(1977)Gibbons, Olkin, and Sobel] page 51 ⟨*square-mean-root (SMR)* (2.4.1)⟩, ✎ [Euclid(circa 300BC)] ⟨triangle inequality—Book I Proposition 20⟩

[39] **metric space:** ✎ [Dieudonné(1969)], page 28, ✎ [Copson(1968)], page 21, ✎ [Hausdorff(1937)] page 109, ✎ [Fréchet(1928)], ✎ [Fréchet(1906)] page 30 **near metric space:** ▣ [Czerwik(1993)] page 5 ⟨*b-metric*; (1),(2),(5)⟩, ▣ [Fagin et al.(2003a)Fagin, Kumar, and Sivakumar], ▣ [Fagin et al.(2003b)Fagin, Kumar, and Sivakumar] ⟨Definition 4.2 (Relaxed metrics)⟩, ▣ [Xia(2009)] page 453 ⟨Definition 2.1.⟩, ✎ [Heinonen(2001)] page 109 ⟨14.1 Quasimetric spaces.⟩, ✎ [Kirk and Shahzad(2014)] page 113 ⟨Definition 12.1⟩, ✎ [Deza and Deza(2014)] page 7

$\le$ *and usual metric* $|\cdot|$ *on* $\mathbb{R}$. *The function* $\tau(p, \sigma; x, y, z; \mathsf{d})$ *is* CONTINUOUS *and* STRICTLY MONOTONE *in* $(\mathbb{R}, |\cdot|, \le)$ *with respect to both the variables* $p$ *and* $\sigma$.

✎PROOF:

(1) Proof that $\tau(p, \sigma; x, y, z; \mathsf{d})$ is *continuous* and *strictly monotone* with respect to $p$: This follows directly from Theorem B.7 (page 31).

(2) Proof that $\tau(p, \sigma; x, y, z; \mathsf{d})$ is *continuous* and *strictly monotone* with respect to $\sigma$:

$$\tau(p, \sigma; x, y, z; \mathsf{d}) \triangleq 2\sigma \underbrace{\left[\frac{1}{2}\mathsf{d}^p(x, z) + \frac{1}{2}\mathsf{d}^p(z, y)\right]^{\frac{1}{p}}}_{\mathsf{f}(p, x, y, z)} \qquad \text{by definition of } \tau \text{ (Definition 4.1 page 15)}$$

$$= 2\sigma\mathsf{f}(p, x, y, z) \qquad \text{where } \mathsf{f} \text{ is defined as above}$$

$\implies \tau$ is *affine* with respect to $\sigma$

$\implies \tau$ is *continuous* and *strictly monotone* with respect to $\sigma$:

☞

**Corollary 4.7** *Let* $\tau(p, \sigma; x, y, z; \mathsf{d})$ *be the* POWER TRIANGLE FUNCTION *in the* DISTANCE SPACE *(Definition 3.1 page 6)* $(X, \mathsf{d})$.

$$\tau(p, \sigma; x, y, z; \mathsf{d}) = \begin{cases} 2\sigma \max\{\mathsf{d}(x, z), \mathsf{d}(z, y)\} & \text{for} \quad p = \infty, & \text{(MAXIMUM, } corresponds\ to\ \text{INFRAMETRIC SPACE)} \\ 2\sigma\left[\frac{1}{2}\mathsf{d}^2(x, z) + \frac{1}{2}\mathsf{d}^2(z, y)\right]^{\frac{1}{2}} & \text{for} \quad p = 2, & \text{(QUADRATIC MEAN)} \\ \sigma[\mathsf{d}(x, z) + \mathsf{d}(z, y)] & \text{for} \quad p = 1, & \text{(ARITHMETIC MEAN, } corresponds\ to\ \text{NEAR METRIC SPACE)} \\ 2\sigma\sqrt{\mathsf{d}(x, z)}\sqrt{\mathsf{d}(z, y)} & \text{for} \quad p = 0, & \text{(GEOMETRIC MEAN)} \\ 4\sigma\left[\frac{1}{\mathsf{d}(x,z)} + \frac{1}{\mathsf{d}(z,y)}\right]^{-1} & \text{for} \quad p = -1, & \text{(HARMONIC MEAN)} \\ 2\sigma \min\{\mathsf{d}(x, z), \mathsf{d}(z, y)\} & \text{for} \quad p = -\infty, & \text{(MINIMUM)} \end{cases}$$

✎PROOF:    These follow directly from Theorem B.7 (page 31).    ☞

**Corollary 4.8** *Let* $(X, \mathsf{d})$ *be a* DISTANCE SPACE.

$$2\sigma \min\{\mathsf{d}(x, z), \mathsf{d}(z, y)\} \quad \le \quad 4\sigma\left[\frac{1}{\mathsf{d}(x,z)} + \frac{1}{\mathsf{d}(z,y)}\right]^{-1} \quad \le \quad 2\sigma\sqrt{\mathsf{d}(x, z)}\sqrt{\mathsf{d}(z, y)}$$

$$\le \quad \sigma[\mathsf{d}(x, z) + \mathsf{d}(z, y)] \quad \le \quad 2\sigma \max\{\mathsf{d}(x, z), \mathsf{d}(z, y)\}$$

✎PROOF:    These follow directly from Corollary B.8 (page 34).    ☞

### 4.2.2  Properties of power distance spaces

The *power triangle inequality* property of a *power distance space* axiomatically endows a metric with an upper bound. Lemma 4.9 (next) demonstrates that there is a complementary lower bound somewhat similar in form to the *power triangle inequality* upper bound. In the special case where $2\sigma = 2^{\frac{1}{p}}$, the lower bound helps provide a simple proof of the *continuity* of a large class of *power distance function*s *(Theorem 4.18 page 23)*. The inequality $2\sigma \le 2^{\frac{1}{p}}$ is a special relation in this paper and appears repeatedly in this paper; it appears as an inequality in Lemma 4.13 (page 20), Corollary 4.12 (page 20) and Corollary 4.14 (page 21), and as an equality in Lemma 4.9 (next) and Theorem 4.18 (page 23). It is plotted in Figure 1 (page 18).

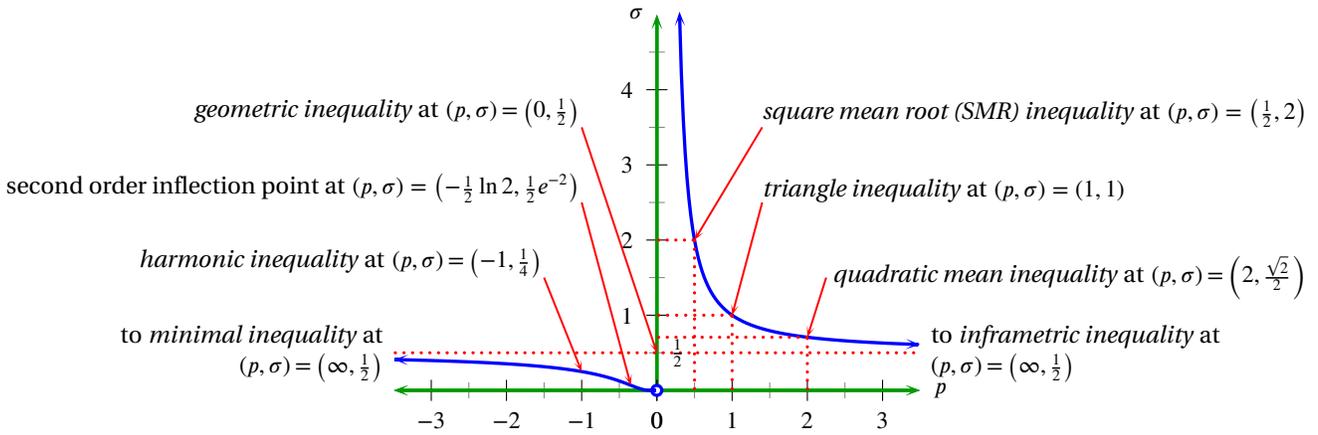

Figure 1: $\sigma = \frac{1}{2}(2^{\frac{1}{p}}) = 2^{\frac{1}{p}-1}$ or $p = \frac{\ln 2}{\ln(2\sigma)}$ (see Lemma 4.9 page 18, Lemma 4.13 page 20, Corollary 4.14 page 21, Corollary 4.12 page 20, and Theorem 4.18 page 23).

**Lemma 4.9** [40] *Let* $(X, \mathsf{d}, p, \sigma)$ *be a* POWER TRIANGLE TRIANGLE SPACE *(Definition 4.3 page 16). Let* $|\cdot|$ *be the* ABSOLUTE VALUE *function (Definition 2.10 page 5). Let* max $\{x, y\}$ *be the maximum and* min $\{x, y\}$ *the minimum of any* $x, y \in \mathbb{R}^*$*. Then, for all* $(p, \sigma) \in \mathbb{R}^* \times \mathbb{R}^+$*,*

1. $\mathsf{d}^p(x, y) \geq \max\left\{0, \frac{2}{(2\sigma)^p}\mathsf{d}^p(x, z) - \mathsf{d}^p(z, x), \frac{2}{(2\sigma)^p}\mathsf{d}^p(y, z) - \mathsf{d}^p(z, x)\right\}$    $\forall x, y, z \in X$    *and*

2. $\mathsf{d}(x, y) \geq |\mathsf{d}(x, z) - \mathsf{d}(z, y)|$    *if* $p \neq 0$ *and* $2\sigma = 2^{\frac{1}{p}}$      $\forall x, y, z \in X$.

✎ PROOF:

(1) lemma: $\frac{2}{(2\sigma)^p}\mathsf{d}^p(x, z) - \mathsf{d}^p(z, y) \leq \mathsf{d}^p(x, y)$    $\forall (p, \sigma) \in \mathbb{R}^* \times \mathbb{R}^+$: Proof:

$$\frac{2}{(2\sigma)^p}\mathsf{d}^p(x, z) - \mathsf{d}^p(z, y) \leq \frac{2}{(2\sigma)^p}\left[2\sigma\left[\tfrac{1}{2}\mathsf{d}^p(x, y) + \tfrac{1}{2}\mathsf{d}^p(y, z)\right]^{\frac{1}{p}}\right]^p - \mathsf{d}^p(z, y) \quad \text{by \textit{power triangle inequality}}$$

$$= \frac{2(2\sigma)^p}{(2\sigma)^p}\left[\tfrac{1}{2}\mathsf{d}^p(x, y) + \tfrac{1}{2}\mathsf{d}^p(y, z)\right] - \mathsf{d}^p(z, y)$$

$$= \left[\mathsf{d}^p(x, y) + \mathsf{d}^p(y, z)\right] - \mathsf{d}^p(y, z) \quad \text{by \textit{symmetric} property of d}$$

$$= \mathsf{d}^p(x, y)$$

(2) Proof for $(p, \sigma) \in \mathbb{R}^* \times \mathbb{R}^+$ case:

$\mathsf{d}^p(x, y) \geq \frac{2}{(2\sigma)^p}\mathsf{d}^p(x, z) - \mathsf{d}^p(z, y)$   by (1) lemma

$\mathsf{d}^p(x, y) = \mathsf{d}^p(y, x) \geq \frac{2}{(2\sigma)^p}\mathsf{d}^p(y, z) - \mathsf{d}^p(z, x)$   by *commutative* property of d and (1) lemma

$\mathsf{d}^p(x, y) \geq 0$   by *non-negative* property of d (Definition 3.1 page 6)

The rest follows because $\mathsf{g}(x) \triangleq x^{\frac{1}{p}}$ is *strictly monotone* in $\mathbb{R}^{\mathbb{R}}$.

---

[40] in *metric space* $((p, \sigma) = (1, 1))$: ✎ [Dieudonné(1969)] page 28, ✎ [Michel and Herget(1993)] page 266, ✎ [Berberian(1961)] page 37 ⟨Theorem II.4.1⟩

(3) Proof for $2\sigma = 2^{\frac{1}{p}}$ case:

$$\mathsf{d}(x,y) \geq \max \left\{ 0,\ \frac{2}{(2\sigma)^p}\mathsf{d}^p(x,z) - \mathsf{d}^p(z,y),\ \frac{2}{(2\sigma)^p}\mathsf{d}^p(y,z) - \mathsf{d}^p(z,x) \right\}^{\frac{1}{p}} \qquad \text{by item (2) (page 18)}$$

$$= \max\{0,\ \mathsf{d}(x,z) - \mathsf{d}(z,y),\ \mathsf{d}(y,z) - \mathsf{d}(z,x)\} \qquad \text{by } 2\sigma = 2^{\frac{1}{p}} \text{ hypothesis} \iff \frac{2}{(2\sigma)^p} = 1$$

$$= \max\{0,\ (\mathsf{d}(x,z) - \mathsf{d}(z,y)),\ -(\mathsf{d}(x,z) - \mathsf{d}(z,y))\} \qquad \text{by } \textit{symmetric} \text{ property of } \mathsf{d}$$

$$= |(\mathsf{d}(x,z) - \mathsf{d}(z,y))|$$

**Theorem 4.10**  *Let* $(X, \mathsf{d}, p, \sigma)$ *be a* POWER DISTANCE SPACE *(Definition 4.3 page 16). Let* $\mathsf{B}$ *be an* OPEN BALL *(Definition 3.5 page 6) on* $(X, \mathsf{d})$. *Then for all* $(p, \sigma) \in (\mathbb{R}^* \setminus \{0\}) \times \mathbb{R}^+$,

$$\left\{ \begin{array}{l} \textit{A.}\ \ 2\sigma \leq 2^{\frac{1}{p}} \quad and \\ \textit{B.}\ \ q \in \mathsf{B}(\theta, r) \end{array} \right\} \implies \left\{ \begin{array}{l} \textit{1.}\ \ \exists r_q \in \mathbb{R}^+ \ \ \text{such that} \\ \quad \mathsf{B}(q, r_q) \subseteq \mathsf{B}(\theta, r) \end{array} \right\} \implies \left\{ \begin{array}{l} \textit{B.}\ \ q \in \mathsf{B}(\theta, r) \end{array} \right\}$$

✎Proof:

(1) lemma:

$$q \in \mathsf{B}(\theta, r) \iff \mathsf{d}(\theta, q) < r \qquad \text{by definition of } \textit{open ball} \text{ (Definition 3.5 page 6)}$$

$$\iff 0 < r - \mathsf{d}(\theta, q) \qquad \text{by field property of real numbers}$$

$$\iff \exists r_q \in \mathbb{R}^+ \ \text{such that}\ \ 0 < r_q < r - \mathsf{d}(\theta, q) \qquad \text{by } \textit{The Archimedean Property}^{[41]}$$

(2) Proof that $(A), (B) \implies (1)$:

$$\mathsf{B}(q, r_q) \triangleq \{x \in X \mid \mathsf{d}(q, x) < r_q\} \qquad \text{by definition of } \textit{open ball} \text{ (Definition 3.5 page 6)}$$

$$= \{x \in X \mid \mathsf{d}^p(q, x) < r_q^p \in \mathbb{R}^+\} \qquad \text{because } \mathsf{f}(x) \triangleq x^p \text{ is } \textit{monotone}$$

$$\subseteq \{x \in X \mid \mathsf{d}^p(q, x) < r^p - \mathsf{d}^p(\theta, q)\} \qquad \text{by hypothesis B and (1) lemma page 19}$$

$$= \{x \in X \mid \mathsf{d}^p(\theta, q) + \mathsf{d}^p(q, x) < r^p\} \qquad \text{by field property of real numbers}$$

$$= \left\{x \in X \mid \left[\mathsf{d}^p(\theta, q) + \mathsf{d}^p(q, x)\right]^{\frac{1}{p}} < r\right\} \qquad \text{because } \mathsf{f}(x) \triangleq x^{\frac{1}{p}} \text{ is } \textit{monotone}$$

$$\subseteq \left\{x \in X \mid 2^{1 - \frac{1}{p}}\sigma \left[\mathsf{d}^p(\theta, q) + \mathsf{d}^p(q, x)\right]^{\frac{1}{p}} < r\right\} \qquad \text{by hypothesis A which implies } 2^{1-\frac{1}{p}}\sigma \leq 1$$

$$= \left\{x \in X \mid 2\sigma \left[\frac{1}{2}\mathsf{d}(\theta, q)\, x + \frac{1}{2}\mathsf{d}^p(q, x)\right]^{\frac{1}{p}} < r\right\} \qquad \text{because } 2^{1-\frac{1}{p}}\sigma = 2\sigma(\frac{1}{2})^{\frac{1}{p}}$$

$$\triangleq \{x \in X \mid \tau(p, \sigma, \theta, x, q) < r\} \qquad \text{by definition of } \tau \text{ (Definition 4.1 page 15)}$$

$$\subseteq \{x \in X \mid \mathsf{d}(\theta, x) < r\} \qquad \text{by definition of } (X, \mathsf{d}, p, \sigma) \text{ (Definition 4.3 page 16)}$$

$$\triangleq \mathsf{B}(\theta, r) \qquad \text{by definition of } \textit{open ball} \text{ (Definition 3.5 page 6)}$$

(3) Proof that $(B) \impliedby (1)$:

$$q \in \{x \in X \mid \mathsf{d}(q, x) = 0\} \qquad \text{by } \textit{nondegenerate} \text{ property (Definition 3.1 page 6)}$$

$$\subseteq \{x \in X \mid \mathsf{d}(q, x) < r_q\} \qquad \text{because } r_q > 0$$

$$\triangleq \mathsf{B}(q, r_q) \qquad \text{by definition of } \textit{open ball} \text{ (Definition 3.5 page 6)}$$

$$\subseteq \mathsf{B}(\theta, r) \qquad \text{by hypothesis 2}$$

---

[41] ✎ [Aliprantis and Burkinshaw(1998)] page 17 ⟨Theorem 3.3 ("*The Archimedean Property*") and Theorem 3.4⟩, [Zorich(2004)] page 53 ⟨6° ("*The principle of Archimedes*") and 7°⟩

**Corollary 4.11**  *Let* $(X, \mathsf{d}, p, \sigma)$ *be a* POWER DISTANCE SPACE. *Then for all* $(p, \sigma) \in (\mathbb{R}^* \setminus \{0\}) \times \mathbb{R}^+$,

$$\left\{ 2\sigma \le 2^{\frac{1}{p}} \right\} \quad \Longrightarrow \quad \left\{ every \text{ OPEN BALL } in \ (X, \mathsf{d}) \ is \text{ OPEN} \right\}$$

✎PROOF:    This follows from Theorem 4.10 (page 19) and Theorem 3.10 (page 8).

**Corollary 4.12**  *Let* $(X, \mathsf{d}, p, \sigma)$ *be a* POWER DISTANCE SPACE. *Let* $\boldsymbol{B}$ *be the set of all* OPEN BALL*s in* $(X, \mathsf{d})$.
*Then for all* $(p, \sigma) \in (\mathbb{R}^* \setminus \{0\}) \times \mathbb{R}^+$,

$$\left\{ 2\sigma \le 2^{\frac{1}{p}} \right\} \quad \Longrightarrow \quad \left\{ \ \boldsymbol{B} \ is \ a \text{ BASE } for \ (X, \boldsymbol{T}) \ \right\}$$

✎PROOF:

(1)  The set of all *open balls* in $(X, \mathsf{d})$ is a *base* for $(X, \boldsymbol{T})$ by Corollary 4.11 (page 20) and Theorem A.4 (page 27).

(2)  $\boldsymbol{T}$ is a topology on $X$ by Definition A.3 (page 27).

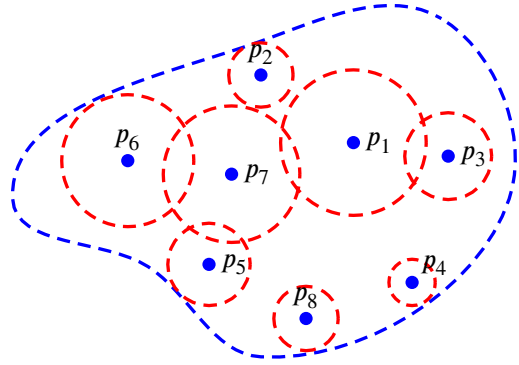

Figure 2: *open set* (see Lemma 4.13 page 20)

Lemma 4.13 (next) demonstrates that every point in an open set is contained in an open ball that is contained in the original open set (see also Figure 2 page 20).

**Lemma 4.13**  *Let* $(X, \mathsf{d}, p, \sigma)$ *be a* POWER DISTANCE SPACE. *Let* $\mathsf{B}$ *be an* OPEN BALL *on* $(X, \mathsf{d})$. *Then for all*
$(p, \sigma) \in (\mathbb{R}^* \setminus \{0\}) \times \mathbb{R}^+$,

$$\left\{ \begin{array}{ll} A. & 2\sigma \le 2^{\frac{1}{p}} \quad and \\ B. & U \ is \text{ OPEN } in \ (X, \mathsf{d}) \end{array} \right\} \Longrightarrow \left\{ \begin{array}{l} I. \ \forall x \in U, \ \exists r \in \mathbb{R}^+ \ \ such \ that \\ \ \ \ \mathsf{B}(x, r) \subseteq U \end{array} \right\} \Longrightarrow \left\{ \begin{array}{l} B. \ U \ is \\ \ \ \ \text{OPEN } in \ (X, \mathsf{d}) \end{array} \right\}$$

✎PROOF:

(1)  Proof that for $((A), (B) \implies (1))$:

$$U = \bigcup \left\{ \mathsf{B}(x_\gamma, r_\gamma) \, \middle| \, \mathsf{B}(x_\gamma, r_\gamma) \subseteq U \right\} \qquad \text{by left hypothesis and Corollary 4.12 page 20}$$

$$\supseteq \mathsf{B}(x, r) \qquad\qquad\qquad\qquad \text{because } x \text{ must be in one of those balls in } U$$

(2) Proof that $((B) \Longleftarrow (1))$ case:

$$U = \bigcup \{ x \in X \,|\, x \in U \} \qquad \text{by definition of union operation } \bigcup$$
$$= \bigcup \big\{ \mathsf{B}(x,r) \,\big|\, x \in U \text{ and } \mathsf{B}(x,r) \subseteq U \big\} \qquad \text{by hypothesis (1)}$$
$$\Longrightarrow U \text{ is } open \qquad \text{by Corollary 4.12 page 20 and Corollary 3.8 page 8}$$

**Corollary 4.14** [42] *Let* $(X, \mathsf{d}, p, \sigma)$ *be a* POWER DISTANCE SPACE. *Let* $\mathsf{B}$ *be an* OPEN BALL *on* $(X, \mathsf{d})$. *Then for all* $(p, \sigma) \in (\mathbb{R}^* \setminus \{0\}) \times \mathbb{R}^+$,

$$\left\{ 2\sigma \le 2^{\frac{1}{p}} \right\} \quad \Longrightarrow \quad \big\{ \text{ every OPEN BALL } \mathsf{B}(x,r) \text{ in } (X, \mathsf{d}) \text{ is OPEN } \big\}$$

✎PROOF:

The union of any set of open balls is open          by Corollary 4.12 page 20

$\qquad \Longrightarrow$ the union of a set of just one open ball is open

$\qquad \Longrightarrow$ every open ball is open.

**Theorem 4.15** [43] *Let* $(X, \mathsf{d}, p, \sigma)$ *be a* POWER DISTANCE SPACE. *Let* $(X, \boldsymbol{T})$ *be a* TOPOLOGICAL SPACE INDUCED BY $(X, \mathsf{d})$. *Let* $(x_n \in X)_{n \in \mathbb{Z}}$ *be a sequence in* $(X, \mathsf{d})$.

$$\underbrace{(x_n) \text{ converges to a limit } x}_{\text{(Definition A.16 page 28)}} \quad \Longleftrightarrow \quad \left\{ \begin{array}{l} \text{for any } \varepsilon \in \mathbb{R}^+, \text{ there exists } N \in \mathbb{Z} \\ \text{such that for all } n > N, \quad \mathsf{d}(x_n, x) < \varepsilon \end{array} \right\}$$

✎PROOF:

$$(x_n) \to x \iff x_n \in U \quad \forall U \in N_x,\, n > N \qquad \text{by Definition A.16 page 28}$$
$$\iff \exists \mathsf{B}(x, \varepsilon) \quad \text{such that} \quad x_n \in \mathsf{B}(x, \varepsilon) \,\forall n > N \qquad \text{by Lemma 4.13 page 20}$$
$$\iff \mathsf{d}(x_n, x) < \varepsilon \qquad \text{by Definition 3.5 page 6}$$

In *distance space*s (Definition 3.1 page 6), not all *convergent* sequences are *Cauchy* (Example 3.22 page 13). However in a distance space with any *power triangle inequality* (Definition 4.3 page 16), all *convergent* sequences are *Cauchy* (next theorem).

**Theorem 4.16** [44] *Let* $(X, \mathsf{d}, p, \sigma)$ *be a* POWER DISTANCE SPACE. *Let* $\mathsf{B}$ *be an* OPEN BALL *on* $(X, \mathsf{d})$. *For any* $(p, \sigma) \in \mathbb{R}^* \times \mathbb{R}^+$,

$$\left\{ \begin{array}{l} (x_n) \text{ is CONVERGENT} \\ \text{in } (X, \mathsf{d}) \end{array} \right\} \Longrightarrow \left\{ \begin{array}{l} (x_n) \text{ is CAUCHY} \\ \text{in } (X, \mathsf{d}) \end{array} \right\} \Longrightarrow \left\{ \begin{array}{l} (x_n) \text{ is BOUNDED} \\ \text{in } (X, \mathsf{d}) \end{array} \right\}$$

---

✎ PROOF:

(1)   Proof that *convergent* $\implies$ *Cauchy*:

$$d(x_n, x_m) \le \tau(p, \sigma; x_n, x_m, x) \qquad \text{by definition of } \textit{power triangle inequality} \text{ (Definition 4.3 page 16)}$$

$$\triangleq 2\sigma\left[\frac{1}{2}d^p(x_n, x) + \frac{1}{2}d^p(x_m, x)\right]^{\frac{1}{p}} \qquad \text{by definition of } \textit{power triangle function} \text{ (Definition 4.1 page 15)}$$

$$< 2\sigma\left[\frac{1}{2}\varepsilon^p + \frac{1}{2}\varepsilon^p\right]^{\frac{1}{p}} \qquad \text{by } \textit{convergence} \text{ hypothesis (Definition A.16 page 28)}$$

$$= 2\sigma\varepsilon \qquad \text{by definition of } \textit{convergence} \text{ (Definition A.16 page 28)}$$

$$\implies \textit{Cauchy} \qquad \text{by definition of } \textit{Cauchy} \text{ (Definition 3.12 page 9)}$$

$$d(x_n, x_m) \le \tau(\infty, \sigma; x_n, x_m, x) \qquad \text{by definition of } \textit{power triangle inequality} \text{ at } p = \infty$$

$$= 2\sigma \max\{d(x_n, x), d(x_m, x)\} \qquad \text{by Corollary 4.7 (page 17)}$$

$$= 2\sigma \max\{\varepsilon, \varepsilon\} \qquad \text{by } \textit{convergent} \text{ hypothesis (Definition A.16 page 28)}$$

$$= 2\sigma\varepsilon \qquad \text{by definition of max}$$

$$d(x_n, x_m) \le \tau(-\infty, \sigma; x_n, x_m, x) \qquad \text{by definition of } \textit{power triangle inequality} \text{ at } p = -\infty$$

$$= 2\sigma \min\{d(x_n, x), d(x_m, x)\} \qquad \text{by Corollary 4.7 (page 17)}$$

$$= 2\sigma \min\{\varepsilon, \varepsilon\} \qquad \text{by } \textit{convergent} \text{ hypothesis (Definition A.16 page 28)}$$

$$= 2\sigma\varepsilon \qquad \text{by definition of min}$$

$$d(x_n, x_m) \le \tau(0, \sigma; x_n, x_m, x) \qquad \text{by definition of } \textit{power triangle inequality} \text{ at } p = 0$$

$$= 2\sigma\sqrt{d(x_n, x)}\sqrt{d(x_m, x)} \qquad \text{by Corollary 4.7 (page 17)}$$

$$= 2\sigma\sqrt{\varepsilon}\sqrt{\varepsilon} \qquad \text{by } \textit{convergent} \text{ hypothesis (Definition A.16 page 28)}$$

$$= 2\sigma\varepsilon \qquad \text{by property of } \mathbb{R}$$

(2)   Proof that *Cauchy* $\implies$ *bounded*: by Proposition 3.14 (page 9).

**Theorem 4.17** [45] Let $(X, d, p, \sigma)$ be a POWER DISTANCE SPACE. Let $f \in \mathbb{Z}^{\mathbb{Z}}$ be a STRICTLY MONOTONE *function such that* $f(n) < f(n+1)$. *For any* $(p, \sigma) \in \mathbb{R}^* \times \mathbb{R}^+$

$$\left\{ \begin{array}{ll} \textit{1.} & (x_n)_{n \in \mathbb{Z}} \textit{ is CAUCHY} \quad \textit{and} \\ \textit{2.} & (x_{f(n)})_{n \in \mathbb{Z}} \textit{ is CONVERGENT} \end{array} \right\} \implies \left\{ (x_n)_{n \in \mathbb{Z}} \textit{ is CONVERGENT.} \right\}$$

✎ PROOF:

$$d(x_n, x) = d(x, x_n) \qquad \text{by } \textit{symmetric} \text{ property of d}$$

$$\le \tau(p, \sigma; x, x_n, x_{f(n)}) \qquad \text{by definition of } \textit{power triangle inequality} \text{ (Definition 4.3 page 16)}$$

$$\triangleq 2\sigma\left[\frac{1}{2}d^p(x, x_{f(n)}) + \frac{1}{2}d^p(x_{f(n)}, x_n)\right]^{\frac{1}{p}} \qquad \text{by definition of } \textit{power triangle function} \text{ (Definition 4.1 page 15)}$$

$$= 2\sigma\left[\frac{1}{2}\varepsilon + \frac{1}{2}d^p(x_{f(n)}, x_n)\right]^{\frac{1}{p}} \qquad \text{by left hypothesis 2}$$

$$= 2\sigma\left[\frac{1}{2}\varepsilon^p + \frac{1}{2}\varepsilon^p\right]^{\frac{1}{p}} \qquad \text{by left hypothesis 1}$$

$$= 2\sigma\varepsilon$$

$$\implies \textit{convergent} \qquad \text{by definition of } \textit{convergent} \text{ (Definition A.16 page 28)}$$

[45] in *metric space*: ✎ [Rosenlicht(1968)] page 52

**Theorem 4.18** [46] *Let* $(X, d, p, \sigma)$ *be a* POWER DISTANCE SPACE. *Let* $(\mathbb{R}, q)$ *be a metric space of real numbers with the usual metric* $q(x, y) \triangleq |x - y|$. *Then*

$$\left\{ 2\sigma = 2^{\frac{1}{p}} \right\} \implies \left\{ d \text{ is CONTINUOUS in } (\mathbb{R}, q) \right\}$$

✎PROOF:

$$
\begin{aligned}
\left| d(x, y) - d(x_n, y_n) \right| &\le \left| d(x, y) - d(x_n, y) \right| + \left| d(x_n, y) - d(x_n, y_n) \right| && \text{by } \textit{triangle inequality} \text{ of } (\mathbb{R}, |x - y|) \\
&= \left| d(x, y) - d(y, x_n) \right| + \left| d(y, x_n) - d(x_n, y_n) \right| && \text{by } \textit{commutative} \text{ property of d (Definition 3.1 page 6)} \\
&\le d(x, x_n) + d(y, y_n) && \text{by } 2\sigma = 2^{\frac{1}{p}} \text{ and Lemma 4.9 (page 18)} \\
&= 0 && \text{as } n \to \infty
\end{aligned}
$$

In *distance space*s and *topological space*s, limits of convergent sequences are in general *not unique* (Example 3.21 page 12, Example A.17 page 29). However Theorem 4.19 (next) demonstrates that, in a *power distance space*, limits *are* unique.

**Theorem 4.19** (Uniqueness of limit) [47] *Let* $(X, d, p, \sigma)$ *be a* POWER DISTANCE SPACE. *Let* $x, y, \in X$ *and let* $(x_n \in X)$ *be an* $X$-*valued sequence.*

$$
\left\{
\begin{array}{ll}
1. & \left\{ \left( (x_n), (y_n) \right) \to (x, y) \right\} \quad \textit{and} \\
2. & (p, \sigma) \in \mathbb{R}^* \times \mathbb{R}^+
\end{array}
\right\}
\implies \left\{ x = y \right\}
$$

✎PROOF:

(1)   lemma: Proof that for all $(p, \sigma) \in \mathbb{R}^* \times \mathbb{R}^+$ and for any $\varepsilon \in \mathbb{R}^+$, there exists $N$ such that $d(x, y) < 2\sigma\varepsilon$:

$$
\begin{aligned}
d(x, y) &\le \tau(p, \sigma; x, y, x_n) && \text{by definition of } \textit{power triangle inequality} \text{ (Definition 4.3 page 16)} \\
&\triangleq 2\sigma \left[ \frac{1}{2} d^p(x, x_n) + \frac{1}{2} d^p(x_n, y) \right]^{\frac{1}{p}} && \text{by definition of } \textit{power triangle function} \text{ (Definition 4.1 page 15)} \\
&< 2\sigma \left[ \frac{1}{2} \varepsilon^p + \frac{1}{2} \varepsilon^p \right]^{\frac{1}{p}} && \text{by left hypothesis and for } p \in \mathbb{R}^* \setminus \{-\infty, 0, \infty\} \\
&= 2\sigma\varepsilon \\
d(x, y) &\le \tau(\infty, \sigma; x, y, x_n) && \text{by definition of } \textit{power triangle inequality} \text{ at } p = \infty \\
&= 2\sigma \max \left\{ d(x, x_n), d(x_n, y) \right\} && \text{by Corollary 4.7 (page 17)} \\
&< 2\sigma\varepsilon && \text{by left hypothesis} \\
d(x, y) &\le \tau(-\infty, \sigma; x, y, x_n) && \text{by definition of } \textit{power triangle inequality} \text{ at } p = -\infty \\
&= 2\sigma \min \left\{ d(x, x_n), d(x_n, y) \right\} && \text{by Corollary 4.7 (page 17)} \\
&< 2\sigma\varepsilon && \text{by left hypothesis} \\
d(x, y) &\le \tau(0, \sigma; x, y, x_n) && \text{by definition of } \textit{power triangle inequality} \text{ at } p = 0 \\
&= 2\sigma \sqrt{d(x, x_n)} \sqrt{d(x_n, y)} && \text{by Corollary 4.7 (page 17)} \\
&= 2\sigma \sqrt{\varepsilon} \sqrt{\varepsilon} && \text{by left hypothesis} \\
&< 2\sigma\varepsilon && \text{by property of real numbers}
\end{aligned}
$$

---

[46] in *metric space* $((p, \sigma) = (1, 1)$ case): ✎ [Berberian(1961)] page 37 ⟨Theorem II.4.1⟩
[47] in *metric space*: ✎ [Rosenlicht(1968)] page 46, ✎ [Thomson et al.(2008)Thomson, Bruckner, and Bruckner] page 32 ⟨Theorem 2.8⟩

(2) Proof that $x = y$ (proof by contradiction):

$$x \neq y \implies \mathsf{d}(x, y) \neq 0 \qquad \text{by the \textit{nondegenerate} property of } \mathsf{d} \text{ (Definition 3.1 page 6)}$$
$$\implies \mathsf{d}(x, y) > 0 \qquad \text{by \textit{non-negative} property of } \mathsf{d} \text{ (Definition 3.1 page 6)}$$
$$\implies \exists \varepsilon \quad \text{such that} \quad \mathsf{d}(x, y) > 2\sigma\varepsilon$$
$$\implies \textit{contradiction} \text{ to (1) lemma page 23}$$
$$\implies \mathsf{d}(x, y) = 0$$
$$\implies x = y$$

## 4.3 Examples

It is not always possible to find a *triangle relation* (Definition 4.3 page 16) $\bigotimes(p, \sigma; \mathsf{d})$ that holds in every *distance space* (Definition 3.1 page 6), as demonstrated by Example 4.20 and Example 4.21 (next two examples).

**Example 4.20** Let $\mathsf{d}(x, y) \in \mathbb{R}^{\mathbb{R} \times \mathbb{R}}$ be defined such that

$$\mathsf{d}(x, y) \triangleq \left\{ \begin{array}{ll} y & \forall (x, y) \in \{4\} \times (0 : 2] \quad \textit{(vertical half-open interval)} \\ x & \forall (x, y) \in (0 : 2] \times \{4\} \quad \textit{(horizontal half-open interval)} \\ |x - y| & \text{otherwise} \quad \textit{(Euclidean)} \end{array} \right\}.$$

Note the following about the pair $(\mathbb{R}, \mathsf{d})$:

(1) By Example 3.21 (page 12), $(\mathbb{R}, \mathsf{d})$ is a *distance space*, but not a *metric space*—that is, the *triangle relation* $\bigotimes(1, 1; \mathsf{d})$ does not hold in $(\mathbb{R}, \mathsf{d})$.

(2) Observe further that $(\mathbb{R}, \mathsf{d})$ is *not a power distance space*. In particular, the *triangle relation* $\bigotimes(p, \sigma; \mathsf{d})$ does not hold in $(\mathbb{R}, \mathsf{d})$ for any finite value of $\sigma$ (does not hold for any $\sigma \in \mathbb{R}^+$):

$$\mathsf{d}(0, 4) = 4 \nleq 0 = \lim_{\varepsilon \to 0} 2\sigma\varepsilon = \lim_{\varepsilon \to 0} 2\sigma \left[ \tfrac{1}{2}|0 - \varepsilon|^p + \tfrac{1}{2}\varepsilon^p \right]^{\frac{1}{p}}$$
$$\triangleq \lim_{\varepsilon \to 0} 2\sigma \left[ \tfrac{1}{2}\mathsf{d}^p(0, \varepsilon) + \tfrac{1}{2}\mathsf{d}^p(\varepsilon, 4) \right]^{\frac{1}{p}} \triangleq \lim_{\varepsilon \to 0} \bigotimes(p, \sigma; 0, 4, \varepsilon; \mathsf{d})$$

**Example 4.21** Let $\mathsf{d}(x, y) \in \mathbb{R}^{\mathbb{R} \times \mathbb{R}}$ be defined such that

$$\mathsf{d}(x, y) \triangleq \left\{ \begin{array}{ll} |x - y| & \text{for } x = 0 \text{ or } y = 0 \text{ or } x = y \quad \textit{(Euclidean)} \\ 1 & \text{otherwise} \quad \textit{(discrete)} \end{array} \right\}.$$

Note the following about the pair $(\mathbb{R}, \mathsf{d})$:

(1) By Example 3.22 (page 13), $(\mathbb{R}, \mathsf{d})$ is a *distance space*, but not a *metric space*—that is, the *triangle relation* $\bigotimes(1, 1; \mathsf{d})$ does not hold in $(\mathbb{R}, \mathsf{d})$.

(2) Observe further that $(\mathbb{R}, \mathsf{d})$ is *not a power distance space*—that is, the *triangle relation* $\bigotimes(p, \sigma; \mathsf{d})$ does not hold in $(\mathbb{R}, \mathsf{d})$ for any value of $(p, \sigma) \in \mathbb{R}^* \times \mathbb{R}^+$:

(a) Proof that $\bigotimes(p, \sigma; \mathsf{d})$ does not hold for any $(p, \sigma) \in \{\infty\} \times \mathbb{R}^+$:

$$\lim_{n,m \to \infty} \mathsf{d}(1/n, 1/m) \triangleq 1 \nleq 0 = 2\sigma \max\{0, 0\} \qquad \text{by definition of } \mathsf{d}$$
$$= 2\sigma \lim_{n,m \to \infty} \max\{\mathsf{d}(1/n, 0), \mathsf{d}(0, 1/m)\} \qquad \text{by Corollary 4.7 (page 17)}$$
$$\geq \lim_{n,m \to \infty} 2\sigma \left[ \tfrac{1}{2}\mathsf{d}^p(1/n, 0) + \tfrac{1}{2}\mathsf{d}^p(0, 1/m) \right]^{\frac{1}{p}} \qquad \text{by Corollary 4.6 (page 16)}$$
$$\triangleq \lim_{n,m \to \infty} \tau(p, \sigma, 1/n, 1/m, 0) \qquad \text{by definition of } \tau \text{ (Definition 4.1 page 15)}$$

(b) Proof that $\oslash(p, \sigma; \mathsf{d})$ does not hold for any $(p, \sigma) \in \mathbb{R}^* \times \mathbb{R}^+$: By Corollary 4.6 (page 16), the *triangle function* (Definition 4.1 page 15) $\tau(p, \sigma; x, y, z; \mathsf{d})$ is *continuous* and *strictly monotone* in $(\mathbb{R}, |\cdot|, \leq)$ with respect to the variable $p$. Item 2a demonstrates that $\oslash(p, \sigma; \mathsf{d})$ fails to hold at the best case of $p = \infty$, and so by Corollary 4.6, it doesn't hold for any other value of $p \in \mathbb{R}^*$ either.

**Example 4.22** Let $\mathsf{d}$ be a function in $\mathbb{R}^{\mathbb{R} \times \mathbb{R}}$ such that

$$\mathsf{d}(x, y) \triangleq \left\{ \begin{array}{ll} 2|x - y| & \forall (x, y) \in \{(0, 1), (1, 0)\} \quad \text{(dilated Euclidean)} \\ |x - y| & \text{otherwise} \quad\quad\quad\quad\quad \text{(Euclidean)} \end{array} \right\}.$$

Note the following about the pair $(\mathbb{R}, \mathsf{d})$:

(1) By Example 3.23 (page 14), $(\mathbb{R}, \mathsf{d})$ is a *distance space*, but not a *metric space*—that is, the *triangle relation* $\oslash(1, 1; \mathsf{d})$ does not hold in $(\mathbb{R}, \mathsf{d})$.

(2) But observe further that $(\mathbb{R}, \mathsf{d}, 1, 2)$ *is a power distance space*:

    (a) Proof that $\oslash(1, 2; \mathsf{d})$ (Definition 4.3 page 16) holds for all $(x, y) \in \{(0, 1), (1, 0)\}$:

$$\begin{aligned}
\mathsf{d}(1, 0) = \mathsf{d}(0, 1) &\triangleq 2|0 - 1| = 2 & \text{by definition of } \mathsf{d} \\
&\leq 2 \leq 2(|0 - z| + |z - 1|) \quad \forall z \in \mathbb{R} & \text{by definition of } |\cdot| \text{ (Definition 2.10 page 5)} \\
&= 2\sigma(\tfrac{1}{2}|0 - z|^p + \tfrac{1}{2}|z - 1|^p)^{\frac{1}{p}} \quad \forall z \in \mathbb{R} & \text{for } (p, \sigma) = (1, 2) \\
&\triangleq 2\sigma(\tfrac{1}{2}\mathsf{d}^p(0, z) + \mathsf{d}^p(z, 1))^{\frac{1}{p}} \quad \forall z \in \mathbb{R} & \text{for } (p, \sigma) = (1, 2) \text{ and by definition of } \mathsf{d} \\
&\triangleq \tau(1, 2; 0, 1, z) & \text{by definition of } \tau \text{ (Definition 4.1 page 15)}
\end{aligned}$$

    (b) Proof that $\oslash(1, 2; \mathsf{d})$ holds for all other $(x, y) \in \mathbb{R}^* \times \mathbb{R}^+$:

$$\begin{aligned}
\mathsf{d}(x, y) &\triangleq 2|x - y| & \text{by definition of } \mathsf{d} \\
&\leq (|x - z| + |z - y|) & \text{by property of } \textit{Euclidean metric space}\text{s} \\
&= 2\sigma(\tfrac{1}{2}|0 - z|^p + \tfrac{1}{2}|z - 1|^p)^{\frac{1}{p}} & \text{for } (p, \sigma) = (1, 1) \\
&\triangleq \tau(1, 1; x, y, z) & \text{by definition of } \tau \text{ (Definition 4.1 page 15)} \\
&\leq \tau(1, 2; x, y, z) & \text{by Corollary 4.6 (page 16)}
\end{aligned}$$

(3) In $(X, \mathsf{d})$, the limits of *convergent* sequences are *unique*. This follows directly from the fact that $(\mathbb{R}, \mathsf{d}, 1, 2)$ is a *power distance space* (item (2) page 25) and by Theorem 4.19 page 23.

(4) In $(X, \mathsf{d})$, *convergent* sequences are *Cauchy*. This follows directly from the fact that $(\mathbb{R}, \mathsf{d}, 1, 2)$ is a *power distance space* (item (2) page 25) and by Theorem 4.16 page 21.

**Example 4.23** Let $\mathsf{d}$ be a function in $\mathbb{R}^{\mathbb{R} \times \mathbb{R}}$ such that $\mathsf{d}(x, y) \triangleq (x - y)^2$. Note the following about the pair $(\mathbb{R}, \mathsf{d})$:

(1) It was demonstrated in Example 3.24 (page 15) that $(\mathbb{R}, \mathsf{d})$ is a *distance space*, but that it is *not* a *metric space* because the *triangle inequality* does not hold.

(2) However, the tuple $(\mathbb{R}, \mathsf{d}, p, \sigma)$ *is a power distance space* (Definition 4.3 page 16) for any $(p, \sigma) \in \mathbb{R}^* \times [2 : \infty)$: In particular, for all $x, y, z \in \mathbb{R}$, the *power triangle inequality* (Definition 4.3 page 16) must hold. The "worst case" for this is when a third point $z$ is exactly "halfway between" $x$ and $y$ in $\mathsf{d}(x, y)$;

that is, when $z = \frac{x+y}{2}$:

$$
\begin{aligned}
(x - y)^2 &\triangleq \mathrm{d}(x, y) && \text{by definition of } \mathrm{d} \\
&\leq \tau(p, \sigma; x, y, z; \mathrm{d}) && \text{by definition } \textit{power triangle inequality} \\
&\triangleq 2\sigma \left[ \tfrac{1}{2}\mathrm{d}^p(x, z) + \tfrac{1}{2}\mathrm{d}^p(z, y) \right]^{\frac{1}{p}} && \text{by definition } \tau \text{ (Definition 4.1 page 15)} \\
&\triangleq 2\sigma \left[ \tfrac{1}{2}(x - z)^{2p} + \tfrac{1}{2}(z - y)^{2p} \right]^{\frac{1}{p}} && \text{by definition of } \mathrm{d} \\
&= 2\sigma \left[ \tfrac{1}{2}|x - z|^{2p} + \tfrac{1}{2}|z - y|^{2p} \right]^{\frac{1}{p}} && \text{because } (x)^2 = |x|^2 \text{ for all } x \in \mathbb{R} \\
&= 2\sigma \left[ \tfrac{1}{2}\left|x - \frac{x+y}{2}\right|^{2p} + \tfrac{1}{2}\left|\frac{x+y}{2} - y\right|^{2p} \right]^{\frac{1}{p}} && \text{because } z = \frac{x+y}{2} \text{ is the "worst case" scenario} \\
&= 2\sigma \left[ \tfrac{1}{2}\left|\frac{y-x}{2}\right|^{2p} + \tfrac{1}{2}\left|\frac{x-y}{2}\right|^{2p} \right]^{\frac{1}{p}} \\
&= 2\sigma \left[ \left|\frac{x-y}{2}\right|^{2p} \right]^{\frac{1}{p}} = \frac{2\sigma}{4}|x - y|^2 \\
&\implies (p, \sigma) \in \mathbb{R}^* \times [2 : \infty)
\end{aligned}
$$

(3) The *power distance function* $\mathrm{d}$ is *continuous* in $(\mathbb{R}, \mathrm{d}, p, \sigma)$ for any $(p, \sigma)$ such that $\sigma \geq 2$ and $2\sigma = p^{\frac{1}{p}}$. This follows directly from Theorem 4.18 (page 23).

# Appendix A   Topological Spaces

**Definition A.1** [48]   Let $\Gamma$ be a set with an arbitrary (possibly uncountable) number of elements. Let $2^X$ be the *power set* of a set $X$ (Definition 2.2 page 4). A family of sets $T \subseteq 2^X$ is a **topology** on $X$ if

1. $\emptyset \in T$ and
2. $X \in T$ and
3. $U, V \in T$ $\implies U \cap V \in T$ and
4. $\{U_\gamma \,|\, \gamma \in \Gamma\} \subseteq T$ $\implies \bigcup_{\gamma \in \Gamma} U_\gamma \in T$ .

The ordered pair $(X, T)$ is a *topological space* if $T$ is a *topology* on $X$. A set $U$ is **open** in $(X, T)$ if $U$ is any element of $T$. A set $D$ is **closed** in $(X, T)$ if $D^c$ is *open* in $(X, T)$.

Just as the power set $2^X$ and the set $\{\emptyset, X\}$ are algebras of sets on a set $X$, so also are these sets topologies on $X$ (next example):

**Example A.2** [49]   Let $\mathcal{T}(X)$ be the set of topologies on a set $X$ and $2^X$ the *power set* (Definition 2.2 page 4) on $X$.

$\{\emptyset, X\}$   is a *topology* in   $\mathcal{T}(X)$   (*indiscrete topology* or *trivial topology*)
$2^X$   is a *topology* in   $\mathcal{T}(X)$   (*discrete topology*)

**Definition A.3** [50] Let $(X, T)$ be a *topological space*. A set $B \subseteq 2^X$ is a **base** for $T$ if

1. $B \subseteq T$ and
2. $\forall U \in T, \quad \exists \{B_\gamma \in B\}$ such that $U = \bigcup_\gamma B_\gamma$

**Theorem A.4** [51] *Let $(X, T)$ be a* TOPOLOGICAL SPACE. *Let $B$ be a subset of $2^X$ such that $B \subseteq 2^X$.*

$$\{ B \text{ is a BASE for } T \} \iff \left\{ \begin{array}{l} \text{For every } x \in X \text{ and for every OPEN SET } U \text{ containing } x, \\ \text{there exists } B_x \in B \text{ such that } \quad x \in B_x \subseteq U. \end{array} \right\}$$

**Theorem A.5** [52] *Let $(X, T)$ be a* TOPOLOGICAL SPACE *(Definition A.1 page 26) and $B \subseteq 2^X$.*

$$B \text{ is a base for } (X, T) \iff \left\{ \begin{array}{ll} 1. & x \in X \implies \exists B_x \in B \text{ such that } x \in B_x \text{ and} \\ 2. & B_1, B_2 \in B \implies B_1 \cap B_2 \in B \end{array} \right.$$

**Example A.6** [53] Let $(X, d)$ be a *metric space*. The set $B \triangleq \{B(x, r) \mid x \in X, r \in \mathbb{N}\}$ (the set of all open balls in $(X, d)$) is a *base* for a topology on $(X, d)$.

**Example A.7** (the standard topology on the real line) [54] The set $B \triangleq \{(a : b) \mid a, b \in \mathbb{R}, a < b\}$ is a *base* for the metric space $(\mathbb{R}, |b - a|)$ (the *usual metric space* on $\mathbb{R}$).

**Definition A.8** [55] Let $(X, T)$ be a *topological space* (Definition A.1 page 26). Let $2^X$ be the *power set* of $X$.

The set $A^-$ is the **closure** of $A \in 2^X$ if $A^- \triangleq \bigcap \{D \in 2^X \mid A \subseteq D \text{ and } D \text{ is } closed\}$.

The set $A^\circ$ is the **interior** of $A \in 2^X$ if $A^\circ \triangleq \bigcup \{U \in 2^X \mid U \subseteq A \text{ and } U \text{ is } open\}$.

A point $x$ is a **closure point** of $A$ if $x \in A^-$.

A point $x$ is an **interior point** of $A$ if $x \in A^\circ$.

A point $x$ is an **accumulation point** of $A$ if $x \in (A \setminus \{x\})^-$

A point $x$ in $A^-$ is a **point of adherence** in $A$ or is **adherent** to $A$ if $x \in A^-$.

**Proposition A.9** [56] *Let $(X, T)$ be a* TOPOLOGICAL SPACE *(Definition A.1 page 26). Let $A^-$ be the* CLOSURE, *$A^\circ$ the* INTERIOR, *and $\partial A$ the* BOUNDARY *of a set $A$. Let $2^X$ be the* POWER SET *of $X$.*

1. $A^-$ *is* CLOSED $\quad \forall A \in 2^X$.
2. $A^\circ$ *is* OPEN $\quad \forall A \in 2^X$.

**Lemma A.10** [57] *Let $A^-$ be the* CLOSURE, *$A^\circ$ the* INTERIOR, *and $\partial A$ the* BOUNDARY *of a set $A$ in a topological space $(X, T)$. Let $2^X$ be the* POWER SET *of $X$.*

1. $A^\circ \subseteq A \subseteq A^-$ $\quad \forall A \in 2^X$.
2. $A = A^\circ \iff A$ *is* OPEN $\quad \forall A \in 2^X$.
3. $A = A^- \iff A$ *is* CLOSED $\quad \forall A \in 2^X$.

**Definition A.11** [58] Let $(X, T_x)$ and $(Y, T_y)$ be *topological space*s (Definition A.1 page 26). Let f be a function in $Y^X$. A function $f \in Y^X$ is **continuous** if $\underbrace{U \in T_y}_{open\ in\ (Y, T_y)} \implies \underbrace{f^{-1}(U) \in T_x}_{open\ in\ (X, T_x)}$ .

A function is **discontinuous** in $(X, T_y)^{(X, T_x)}$ if it is not *continuous* in $(X, T_y)^{(X, T_x)}$.

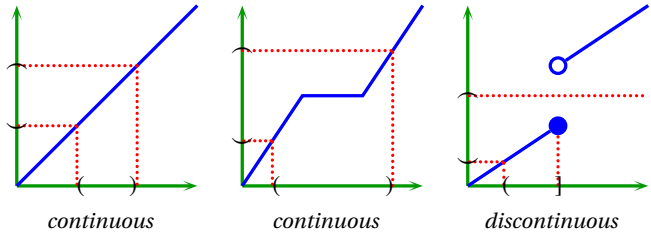

*continuous*          *continuous*          *discontinuous*

Figure 3: *continuous / discontinuous* functions (Example A.12 page 28)

**Example A.12**   Some *continuous / discontinuous* functions are illustrated in Figure 3 (page 28).

Definition A.11 (previous definition) defines continuity using open sets. Continuity can alternatively be defined using closed sets or closure (next theorem).

**Theorem A.13** [59] *Let $(X, T)$ and $(Y, S)$ be topological spaces. Let f be a function in $Y^X$.*
   *The following are equivalent:*
   1. *f is* CONTINUOUS                                            $\iff$
   2. *B is closed in $(Y, S)$ $\implies$ $f^{-1}(B)$ is closed in $(X, T)$*   $\forall B \in 2^Y$  $\iff$
   3. $f(A^-) \subseteq f(A)^-$                               $\forall A \in 2^X$  $\iff$
   4. $f^{-1}(B)^- \subseteq f^{-1}(B^-)$                          $\forall B \in 2^Y$

**Remark A.14**   A word of warning about defining *continuity* in terms of topological spaces—*continuity* is defined in terms of a pair of *topological space*s, and whether function is *continuous* or *discontinuous* in general depends very heavily on the selection of these spaces. This is illustrated in Proposition A.15 (next). The ramification of this is that when declaring a function to be *continuous* or *discontinuous*, one must make clear the assumed *topological space*s.

**Proposition A.15** [60] *Let $(X, T)$ and $(Y, S)$ be* TOPOLOGICAL SPACES*. Let f be a* FUNCTION *in $(Y, S)^{(X, T)}$.*
   1. $T$ *is the* DISCRETE TOPOLOGY   $\implies$   f *is* CONTINUOUS   $\forall f \in (Y, S)^{(X, T)}$
   2. $S$ *is the* INDISCRETE TOPOLOGY   $\implies$   f *is* CONTINUOUS   $\forall f \in (Y, S)^{(X, T)}$

**Definition A.16** [61] Let $(X, T)$ be a *topological space* (Definition A.1 page 26). A sequence $(x_n)_{n \in \mathbb{Z}}$ **converges** in $(X, T)$ to a point $x$ if for each *open set* (Definition A.1 page 26) $U \in T$ that contains $x$ there exists $N \in \mathbb{N}$ such that

$x_n \in U$ for all $n > N$.

This condition can be expressed in any of the following forms:

1. The **limit** of the sequence $(x_n)$ is $x$.
2. The sequence $(x_n)$ is **convergent** with limit $x$.
3. $\lim_{n \to \infty} (x_n) = x$.
4. $(x_n) \to x$.

A sequence that converges is **convergent**. A sequence that does not converge is said to **diverge**, or is **divergent**. An element $x \in A$ is a **limit point** of $A$ if it is the limit of some $A$-valued sequence $(x_n \in A)$.

**Example A.17** [62] Let $(X, T_{31})$ be a *topological space* where $X \triangleq \{x, y, z\}$ and
$$T_{31} \triangleq \{\varnothing, \{x\}, \{x, y\}, \{x, z\}, \{x, y, z\}\}.$$
In this space, the sequence $(x, x, x, \ldots)$ converges to $x$. But this sequence also converges to both $y$ and $z$ because $x$ is in every *open set* (Definition A.1 page 26) that contains $y$ and $x$ is in every *open set* that contains $z$. So, the *limit* (Definition A.16 page 28) of the sequence is *not unique*.

**Example A.18**  In contrast to the low resolution topological space of Example A.17, the limit of the sequence $(x, x, x, \ldots)$ is unique in a *topological space* with sufficiently high resolution with respect to $y$ and $z$ such as the following: Define a *topological space* $(X, T_{56})$ where $X \triangleq \{x, y, z\}$ and
$$T_{56} \triangleq \{\varnothing, \{y\}, \{z\}, \{x, y\}, \{y, z\}, \{x, y, z\}\}.$$
In this space, the sequence $(x, x, x, \ldots)$ converges to $x$ only. The sequence does *not* converge to $y$ or $z$ because there are *open set*s (Definition A.1 page 26) containing $y$ or $z$ that do not contain $x$ (the open sets $\{y\}$, $\{z\}$, and $\{y, z\}$).

**Theorem A.19**  (The Closed Set Theorem) [63] *Let* $(X, T)$ *be a* topological space. *Let* $A$ *be a subset of* $X$ $(A \subseteq X)$. *Let* $A^-$ *be the* closure (Definition A.8 page 27) *of* $A$ *in* $(X, T)$.

$$\underbrace{A \text{ is } \text{CLOSED } in (X, T)}_{(A = A^-)} \quad \Longleftrightarrow \quad \left\{ \begin{array}{l} \textit{Every } A\textit{-valued sequence } (x_n \in A)_{n \in \mathbb{Z}} \\ \textit{that } \text{CONVERGES } in (X, T) \textit{ has its } \text{LIMIT } in A \end{array} \right\}$$

**Theorem A.20** [64] *Let* $(X, T)$ *and* $(Y, S)$ *be* topological spaces. *Let* f *be a function in* $(Y, S)^{(X,T)}$.

$$\underbrace{\left\{ \begin{array}{c} \text{f } \textit{is } \text{CONTINUOUS } in (Y, S)^{(X,T)} \\ \textit{(Definition A.11 page 28)} \end{array} \right\}}_{\text{INVERSE IMAGE CHARACTERIZATION OF CONTINUITY}} \quad \Longleftrightarrow \quad \underbrace{\left\{ \begin{array}{c} (x_n) \to x \implies \text{f}((x_n)) \to \text{f}(x) \\ \textit{(Definition A.16 page 28)} \end{array} \right\}}_{\text{SEQUENTIAL CHARACTERIZATION OF CONTINUITY}}$$

✎Proof:

(1) Proof for the $\implies$ case (proof by contradiction):

(a) Let $U$ be an *open set* in $(Y, T)$ that contains $\text{f}(x)$ but for which there exists no $N$ such that $\text{f}(x_n) \in U$ for all $n > N$.

(b) Note that the set $\text{f}^{-1}(U)$ is also *open* by the *continuity* hypothesis.

(c) If $((x_n)) \to x$, then

$\quad$ $f(((x_n))) \not\to f(x)$ $\implies$ there exists no $N$ such that $f(x_n) \in U$ for all $n > N$ $\qquad$ by Definition A.16 (page 28)

$\qquad\qquad\qquad\quad \implies$ there exists no $M$ such that $x_n \in f^{-1}(U)$ for all $n > M$ $\quad$ by definition of $f^{-1}$

$\qquad\qquad\qquad\quad \implies$ $((x_n)) \not\to x$ $\qquad$ by *continuity* hypothesis and def. of *convergence* (Definition A.16 page 28)

$\qquad\qquad\qquad\quad \implies$ contradiction of $((x_n)) \to x$ hypothesis

$\qquad\qquad\qquad\quad \implies$ $f(((x_n))) \to f(x)$

(2) Proof for the $\impliedby$ case (proof by contradiction):

$\quad$ (a) Let $D$ be a *closed* set in $(Y, \boldsymbol{S})$.

$\quad$ (b) Suppose $f^{-1}(D)$ is *not closed*…

$\quad$ (c) then by the *closed set theorem* (Theorem A.19 page 29), there must exist a *convergent* sequence $((x_n))$ in $(X, \boldsymbol{T})$, but with limit $x$ *not* in $f^{-1}(D)$.

$\quad$ (d) Note that $f(x)$ must be in $D$. Proof:

$\qquad$ (i) by definition of $D$ and $f$, $f(((x_n)))$ is in $D$

$\qquad$ (ii) by left hypothesis, the sequence $f(((x_n)))$ is *convergent* with limit $f(x)$

$\qquad$ (iii) by *closed set theorem* (Theorem A.19 page 29), $f(x)$ must be in $D$.

$\quad$ (e) Because $f(x) \in D$, it must be true that $x \in f^{-1}(D)$.

$\quad$ (f) But this is a contradiction to item (2c) (page 30), and so item (2b) (page 30) must be wrong, and $f^{-1}(D)$ must be *closed*.

$\quad$ (g) And so by Theorem A.13 (page 28), $f$ is *continuous*.

# Appendix B   Finite sums

## B.1   Convexity

**Definition B.1** [65]
A function $f \in \mathbb{R}^{\mathbb{R}}$ is **convex** if

$\quad$ $f\big(\lambda x + [1 - \lambda]y\big) \quad \le \quad \lambda f(x) + (1 - \lambda)\, f(y) \qquad \forall x, y \in \mathbb{R}$ and $\forall \lambda \in (0 : 1)$

A function $g \in \mathbb{R}^{\mathbb{R}}$ is **strictly convex** if

$\quad$ $g\big(\lambda x + [1 - \lambda]y\big) \quad = \quad \lambda g(x) + (1 - \lambda)\, g(y) \qquad \forall x, y \in D, \; x \neq y, \text{ and } \forall \lambda \in (0 : 1)$

A function $f \in \mathbb{R}^{\mathbb{R}}$ is **concave** if $-f$ is *convex*.

A function $f \in \mathbb{R}^{\mathbb{R}}$ is **affine** if $f$ is *convex* and *concave*.

**Theorem B.2** (Jensen's Inequality) [66] *Let* $f \in \mathbb{R}^{\mathbb{R}}$ *be a function.*

$$\left\{ \begin{array}{ll} 1. & f \text{ is CONVEX} \quad \text{\scriptsize(Definition B.1 page 30)} \quad and \\ 2. & \displaystyle\sum_{n=1}^{N} \lambda_n = 1 \qquad \text{\scriptsize(WEIGHTS)} \end{array} \right\} \implies \left\{ f\!\left( \sum_{n=1}^{N} \lambda_n x_n \right) \le \sum_{n=1}^{N} \lambda_n f(x_n) \qquad \forall x_n \in D, \; N \in \mathbb{N} \right\}$$

---

## B.2 Power means

**Definition B.3** [67]

The $\langle \lambda_n \rangle_1^N$ weighted $\phi$-**mean** of a tuple $\langle x_n \rangle_1^N$ is defined as

$$\mathsf{M}_\phi(\langle x_n \rangle) \triangleq \phi^{-1}\left( \sum_{n=1}^N \lambda_n \phi(x_n) \right)$$

where $\phi$ is a *continuous* and *strictly monotonic* function in $\mathbb{R}^{\mathbb{R}^{\vdash}}$

and $\langle \lambda_n \rangle_{n=1}^N$ is a sequence of weights for which $\displaystyle\sum_{n=1}^N \lambda_n = 1$.

**Lemma B.4** [68] *Let* $\mathsf{M}_\phi(\langle x_n \rangle)$ *be the* $\langle \lambda_n \rangle_1^N$ *weighted* $\phi$-*mean and* $\mathsf{M}_\psi(\langle x_n \rangle)$ *the* $\langle \lambda_n \rangle_1^N$ *weighted* $\psi$-*mean of a tuple* $\langle x_n \rangle_1^N$.

$\phi\psi^{-1}$ *is* CONVEX   *and*   $\phi$ *is* INCREASING   $\implies$   $\mathsf{M}_\phi(\langle x_n \rangle) \geq \mathsf{M}_\psi(\langle x_n \rangle)$
$\phi\psi^{-1}$ *is* CONVEX   *and*   $\phi$ *is* DECREASING   $\implies$   $\mathsf{M}_\phi(\langle x_n \rangle) \leq \mathsf{M}_\psi(\langle x_n \rangle)$
$\phi\psi^{-1}$ *is* CONCAVE   *and*   $\phi$ *is* INCREASING   $\implies$   $\mathsf{M}_\phi(\langle x_n \rangle) \leq \mathsf{M}_\psi(\langle x_n \rangle)$
$\phi\psi^{-1}$ *is* CONCAVE   *and*   $\phi$ *is* DECREASING   $\implies$   $\mathsf{M}_\phi(\langle x_n \rangle) \geq \mathsf{M}_\psi(\langle x_n \rangle)$

One of the most well known inequalities in mathematics is *Minkowski's Inequality*. In 1946, H.P. Mulholland submitted a result that generalizes Minkowski's Inequality to an equal weighted $\phi$-mean.[69] And Milovanović and Milovanović (1979) generalized this even further to a *weighted* $\phi$-mean (next).

**Theorem B.5** [70] *Let* $\phi$ *be a function in* $\mathbb{R}^{\mathbb{R}}$.

$$\left\{ \begin{array}{llll} 1. & \phi \text{ is CONVEX} & \text{and} & 2. \quad \phi \text{ is STRICTLY MONOTONE} \quad \text{and} \\ 3. & \phi(0) = 0 & \text{and} & 4. \quad \log \circ \phi \circ \exp \text{ is CONVEX} \end{array} \right\}$$

$$\implies \left\{ \phi^{-1}\left( \sum_{n=1}^N \lambda_n \phi(x_n + y_n) \right) \leq \phi^{-1}\left( \sum_{n=1}^N \lambda_n \phi(x_n) \right) + \phi^{-1}\left( \sum_{n=1}^N \lambda_n \phi(y_n) \right) \right\}$$

**Definition B.6** [71] Let $\mathsf{M}_{\phi(x;p)}(\langle x_n \rangle)$ be the $\langle \lambda_n \rangle_1^N$ weighted $\phi$-mean of a *non-negative* tuple $\langle x_n \rangle_1^N$. A mean $\mathsf{M}_{\phi(x;p)}(\langle x_n \rangle)$ is a **power mean** with parameter $p$ if $\phi(x) \triangleq x^p$. That is,

$$\mathsf{M}_{\phi(x;p)}(\langle x_n \rangle) = \left( \sum_{n=1}^N \lambda_n (x_n)^p \right)^{\frac{1}{p}}$$

**Theorem B.7** [72] *Let* $\mathsf{M}_{\phi(x;p)}(\langle x_n \rangle)$ *be the* POWER MEAN *with parameter* $p$ *of an* $N$-*tuple* $\langle x_n \rangle_1^N$ *in which the elements are* NOT *all equal.*

$$\mathrm{M}_{\phi(x;p)}(\langle x_n\rangle) \;\triangleq\; \left(\sum_{n=1}^{N}\lambda_n(x_n)^p\right)^{\frac{1}{p}} \; is\;\textsc{continuous}\;and\;\textsc{strictly monotone}\;in\;\mathbb{R}^*.$$

$$\mathrm{M}_{\phi(x;p)}(\langle x_n\rangle) \;=\; \begin{cases} \displaystyle\max_{n=1,2,\dots,N}\langle x_n\rangle & for\; p=+\infty \\[2mm] \displaystyle\prod_{n=1}^{N}x_n^{\lambda_n} & for\; p=0 \\[2mm] \displaystyle\min_{n=1,2,\dots,N}\langle x_n\rangle & for\; p=-\infty \end{cases}$$

✎ Proof:

(1) Proof that $M_{\phi(x;p)}$ is *strictly monotone* in $p$:

   (a) Let $p$ and $s$ be such that $-\infty < p < s < \infty$.

   (b) Let $\phi_p \triangleq x^p$ and $\phi_s \triangleq x^s$. Then $\phi_p\phi_s^{-1} = x^{\frac{p}{s}}$.

   (c) The composite function $\phi_p\phi_s^{-1}$ is *convex* or *concave* depending on the values of $p$ and $s$:

| | $p<0$ ($\phi_p$ decreasing) | $p>0$ ($\phi_p$ increasing) |
|---|---|---|
| $s<0$ | *convex* | (not possible) |
| $s>0$ | *convex* | *concave* |

   (d) Therefore by Lemma B.4 (page 31),
$$-\infty < p < s < \infty \implies \mathrm{M}_{\phi(x;p)}(\langle x_n\rangle) < \mathrm{M}_{\phi(x;s)}(\langle x_n\rangle).$$

(2) Proof that $\mathrm{M}_{\phi(x;p)}$ is continuous in $p$ for $p \in \mathbb{R}\setminus 0$: The sum of continuous functions is continuous. For the cases of $p \in \{-\infty,\,0,\,\infty\}$, see the items that follow.

(3) Lemma: $\mathrm{M}_{\phi(x;-p)}(\langle x_n\rangle) = \{\mathrm{M}_{\phi(x;p)}(\langle x_n^{-1}\rangle)\}^{-1}$. Proof:

$$\{\mathrm{M}_{\phi(x;p)}(\langle x_n^{-1}\rangle)\}^{-1} = \left\{\left(\sum_{n=1}^{N}\lambda_n(x_n^{-1})^p\right)^{\frac{1}{p}}\right\}^{-1} \qquad \text{by definition of } \mathrm{M}_\phi$$

$$= \left(\sum_{n=1}^{N}\lambda_n(x_n)^{-p}\right)^{\frac{1}{-p}}$$

$$= \mathrm{M}_{\phi(x;-p)}(\langle x_n\rangle) \qquad \text{by definition of } \mathrm{M}_\phi$$

(4) Proof that $\displaystyle\lim_{p\to\infty}\mathrm{M}_\phi(\langle x_n\rangle) = \max_{n\in\mathbb{Z}}\langle x_n\rangle$:

   (a) Let $x_m \triangleq \displaystyle\max_{n\in\mathbb{Z}}\langle x_n\rangle$

   (b) Note that $\displaystyle\lim_{p\to\infty}\mathrm{M}_\phi \le \max_{n\in\mathbb{Z}}\langle x_n\rangle$ because

$$\lim_{p\to\infty}\mathrm{M}_\phi(\langle x_n\rangle) = \lim_{p\to\infty}\left(\sum_{n=1}^{N}\lambda_n x_n^p\right)^{\frac{1}{p}} \qquad \text{by definition of } \mathrm{M}_\phi$$

$$\le \lim_{p\to\infty}\left(\sum_{n=1}^{N}\lambda_n x_m^p\right)^{\frac{1}{p}} \qquad \begin{array}{l}\text{by definition of } x_m \text{ in item (4a) and because } \phi(x)\triangleq \\ x^p \text{ and } \phi^{-1} \text{ are both increasing or both decreasing}\end{array}$$

$$= \lim_{p\to\infty}\left(x_m^p\underbrace{\sum_{n=1}^{N}\lambda_n}_{1}\right)^{\frac{1}{p}} \qquad \text{because } x_m \text{ is a constant}$$

$$= \lim_{p \to \infty} \left( x_m^p \cdot 1 \right)^{\frac{1}{p}}$$

$$= x_m$$

$$= \max_{n \in \mathbb{Z}} \left( x_n \right) \qquad \text{by definition of } x_m \text{ in item (4a)}$$

(c) But also note that $\lim_{p \to \infty} \mathsf{M}_\phi \geq \max_{n \in \mathbb{Z}} \left( x_n \right)$ because

$$\lim_{p \to \infty} \mathsf{M}_\phi \left( \left( x_n \right) \right) = \lim_{p \to \infty} \left( \sum_{n=1}^{N} \lambda_n x_n^p \right)^{\frac{1}{p}} \qquad \text{by definition of } \mathsf{M}_\phi$$

$$\geq \lim_{p \to \infty} \left( w_m x_m^p \right)^{\frac{1}{p}} \qquad \text{by definition of } x_m \text{ in item (4a) and because } \phi(x) \triangleq x^p \text{ and } \phi^{-1} \text{ are both increasing or both decreasing}$$

$$= \lim_{p \to \infty} w_m^{\frac{1}{p}} x_m^{\frac{p}{p}}$$

$$= x_m$$

$$= \max_{n \in \mathbb{Z}} \left( x_n \right) \qquad \text{by definition of } x_m \text{ in item (4a)}$$

(d) Combining items (b) and (c) we have $\lim_{p \to \infty} \mathsf{M}_\phi = \max_{n \in \mathbb{Z}} \left( x_n \right)$.

(5) Proof that $\lim_{p \to -\infty} \mathsf{M}_\phi \left( \left( x_n \right) \right) = \min_{n \in \mathbb{Z}} \left( x_n \right)$:

$$\lim_{p \to -\infty} \mathsf{M}_{\phi(x;p)} \left( \left( x_n \right) \right) = \lim_{p \to \infty} \mathsf{M}_{\phi(x;-p)} \left( \left( x_n \right) \right) \qquad \text{by change of variable } p$$

$$= \lim_{p \to \infty} \left\{ \mathsf{M}_{\phi(x;p)} \left( \left( x_n^{-1} \right) \right) \right\}^{-1} \qquad \text{by Lemma in item (3) page 32}$$

$$= \lim_{p \to \infty} \frac{1}{\mathsf{M}_{\phi(x;p)} \left( \left( x_n^{-1} \right) \right)}$$

$$= \frac{\lim_{p \to \infty} 1}{\lim_{p \to \infty} \mathsf{M}_{\phi(x;p)} \left( \left( x_n^{-1} \right) \right)} \qquad \text{by property of } \lim \text{ [73]}$$

$$= \frac{1}{\max_{n \in \mathbb{Z}} \left( x_n^{-1} \right)} \qquad \text{by item (4)}$$

$$= \frac{1}{\left( \min_{n \in \mathbb{Z}} \left( x_n \right) \right)^{-1}}$$

$$= \min_{n \in \mathbb{Z}} \left( x_n \right)$$

(6) Proof that $\lim_{p \to 0} \mathsf{M}_\phi \left( \left( x_n \right) \right) = \prod_{n=1}^{N} x_n^{\lambda_n}$:

$$\lim_{p \to 0} \mathsf{M}_\phi \left( \left( x_n \right) \right) = \lim_{p \to 0} \exp \left\{ \ln \left\{ \mathsf{M}_\phi \left( \left( x_n \right) \right) \right\} \right\}$$

$$= \lim_{p \to 0} \exp \left\{ \ln \left\{ \left( \sum_{n=1}^{N} \lambda_n \left( x_n^p \right) \right)^{\frac{1}{p}} \right\} \right\} \qquad \text{by definition of } \mathsf{M}_\phi$$

---

[73]  [Rudin(1976)] page 85 ⟨4.4 Theorem⟩

$$= \exp \left\{ \frac{\frac{\partial}{\partial p} \ln \left( \sum_{n=1}^{N} \lambda_n \left( x_n^p \right) \right)}{\frac{\partial}{\partial p} p} \right\}_{p=0} \qquad \text{by l'Hôpital's rule}^{74}$$

$$= \exp \left\{ \frac{\sum_{n=1}^{N} \lambda_n \frac{\partial}{\partial p} \left( x_n^p \right)}{\sum_{n=1}^{N} \lambda_n \left( x_n^p \right)} \right\}_{p=0} \qquad = \exp \left\{ \frac{\sum_{n=1}^{N} \lambda_n \frac{\partial}{\partial p} \exp \left( \ln \left( x_n^p \right) \right)}{\sum_{n=1}^{N} \lambda_n} \right\}_{p=0}$$

$$= \exp \left\{ \frac{\sum_{n=1}^{N} \lambda_n \frac{\partial}{\partial p} \exp \left( r \ln \left( x_n \right) \right)}{1} \right\}_{p=0} \qquad = \exp \left\{ \sum_{n=1}^{N} \lambda_n \frac{\partial}{\partial p} \exp \left( p \ln \left( x_n \right) \right) \right\}_{p=0}$$

$$= \exp \left\{ \sum_{n=1}^{N} \lambda_n \exp \left\{ p \ln x_n \right\} \ln \left( x_n \right) \right\}_{p=0} \qquad = \exp \left\{ \sum_{n=1}^{N} \lambda_n \ln \left( x_n \right) \right\}$$

$$= \exp \left\{ \sum_{n=1}^{N} \ln \left( x_n^{\lambda_n} \right) \right\} \qquad = \exp \left\{ \ln \prod_{n=1}^{N} x_n^{\lambda_n} \right\} = \prod_{n=1}^{N} x_n^{\lambda_n}$$

⧫

**Corollary B.8** [75]  *Let $\langle\!\langle x_n \rangle\!\rangle_1^N$ be a tuple. Let $\langle\!\langle \lambda_n \rangle\!\rangle_1^N$ be a tuple of weighting values such that $\sum_{n=1}^{N} \lambda_n = 1$.*

$$\min \langle\!\langle x_n \rangle\!\rangle \leq \underbrace{\left( \sum_{n=1}^{N} \lambda_n \frac{1}{x_n} \right)^{-1}}_{\text{harmonic mean}} \leq \underbrace{\prod_{n=1}^{N} x_n^{\lambda_n}}_{\text{geometric mean}} \leq \underbrace{\sum_{n=1}^{N} \lambda_n x_n}_{\text{arithmetic mean}} \leq \max \langle\!\langle x_n \rangle\!\rangle$$

✎ PROOF:

(1) These five means are all special cases of the *power mean* $\mathsf{M}_{\phi(x:p)}$ (Definition B.6 page 31):
  $p = \infty$:    $\max \langle\!\langle x_n \rangle\!\rangle$
  $p = 1$:     arithmetic mean
  $p = 0$:     geometric mean
  $p = -1$:    harmonic mean
  $p = -\infty$:   $\min \langle\!\langle x_n \rangle\!\rangle$

(2) The inequalities follow directly from Theorem B.7 (page 31).

---

(3) Generalized AM-GM inequality: If one is only concerned with the arithmetic mean and geometric mean, their relationship can be established directly using *Jensen's Inequality*:

$$\sum_{n=1}^{N} \lambda_n x_n = b^{\log_b\left(\sum_{n=1}^{N} \lambda_n x_n\right)}$$

$$\geq b^{\left(\sum_{n=1}^{N} \lambda_n \log_b x_n\right)} \qquad \text{by } \textit{Jensen's Inequality} \text{ (Theorem B.2 page 30)}$$

$$= \prod_{n=1}^{N} b^{(\lambda_n \log_b x_n)} = \prod_{n=1}^{N} b^{(\log_b x_n)\lambda_n} = \prod_{n=1}^{N} x_n^{\lambda_n}$$

## B.3   Inequalities

**Lemma B.9** (Young's Inequality)  [76]

$$xy \;<\; \frac{x^p}{p} + \frac{y^q}{q} \quad \text{with } \tfrac{1}{p} + \tfrac{1}{q} = 1 \quad \forall 1 < p < \infty, \; x, y \geq 0, \quad \text{but } y \neq x^{p-1}$$

$$xy \;=\; \frac{x^p}{p} + \frac{y^q}{q} \quad \text{with } \tfrac{1}{p} + \tfrac{1}{q} = 1 \quad \forall 1 < p < \infty, \; x, y \geq 0, \quad \text{and } y = x^{p-1}$$

**Theorem B.10** (Minkowski's Inequality for sequences)  [77] *Let* $(\!(x_n \in \mathbb{C})\!)_1^N$ *and* $(\!(y_n \in \mathbb{C})\!)_1^N$ *be complex* $N$*-tuples.*

$$\left(\sum_{n=1}^{N} |x_n + y_n|^p\right)^{\frac{1}{p}} \leq \left(\sum_{n=1}^{N} |x_n|^p\right)^{\frac{1}{p}} + \left(\sum_{n=1}^{N} |y_n|^p\right)^{\frac{1}{p}} \qquad \forall 1 < p < \infty$$

# Appendix C   Metric preserving functions

**Definition C.1**  [78] Let $\mathbb{M}$ be the set of all *metric space*s (Definition 4.5 page 16) on a set $X$.   $\phi \in \mathbb{R}^{\vdash \mathbb{R}^{\vdash}}$ is a **metric preserving function** if   $\mathrm{d}(x, y) \triangleq \phi \circ \mathrm{p}(x, y)$   is a *metric* on $X$ for all $(X, \mathrm{p}) \in \mathbb{M}$.

**Theorem C.2**  (necessary conditions)  [79] *Let* $\mathcal{R}\phi$ *be the* RANGE *of a function* $\phi$.

$$\left\{\begin{array}{l} \phi \text{ is a} \\ \text{METRIC PRESERVING FUNCTION} \\ \text{(Definition C.1 page 35)} \end{array}\right\} \implies \left\{\begin{array}{ll} 1. & \phi^{-1}(0) = \{0\} & and \\ 2. & \mathcal{R}\phi \subseteq \mathbb{R}^{\vdash} & and \\ 3. & \phi(x + y) \leq \phi(x) + \phi(y) & (\phi \text{ is SUBADDITIVE}) \end{array}\right\}$$

---

**Theorem C.3** (sufficient conditions) [80] *Let $\phi$ be a function in $\mathbb{R}^{\mathbb{R}}$.*

$$\left\{ \begin{array}{ll} 1. & x \geq y \implies \phi(x) \geq \phi(y) \quad \forall x,y \in \mathbb{R}^+ \quad \text{(ISOTONE)} \quad and \\ 2. & \phi(0) = 0 \qquad\qquad\qquad\qquad\qquad\qquad\qquad\qquad\quad and \\ 3. & \phi(x+y) \leq \phi(x) + \phi(y) \quad \forall x,y \in \mathbb{R}^+ \quad \text{(SUBADDITIVE)} \end{array} \right\} \implies \left\{ \begin{array}{l} \phi \text{ is a METRIC} \\ \text{PRESERVING FUNCTION} \\ \textit{(Definition C.1 page 35)}. \end{array} \right\}$$

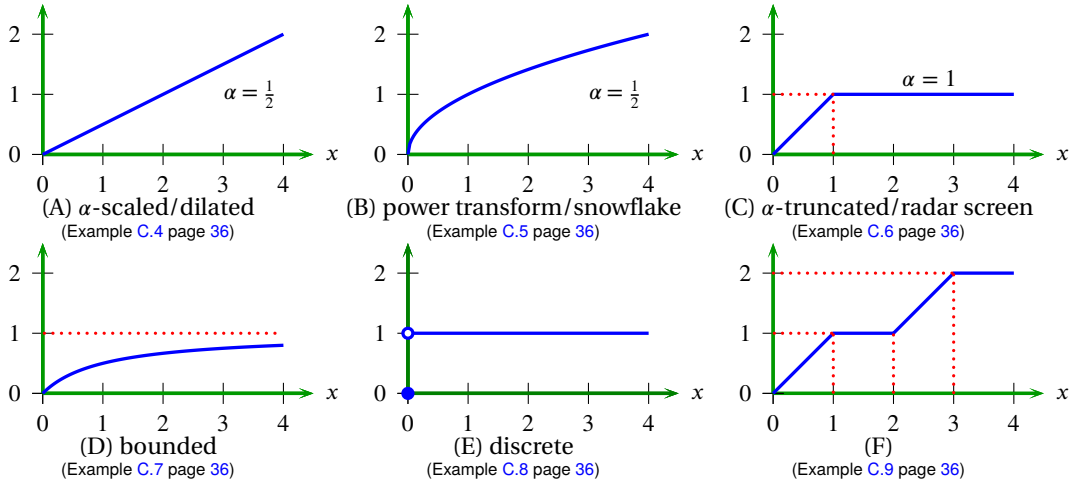

Figure 4: metric preserving functions

**Example C.4** ($\alpha$-scaled metric/dilated metric) [81] Let $(X, \mathsf{d})$ be a *metric space* (Definition 4.5 page 16).
$\phi(x) \triangleq \alpha x, \; \alpha \in \mathbb{R}^+$ is a *metric preserving function* (Figure 4 page 36 (A))

✎ PROOF: The proofs for Example C.4–Example C.9 (page 36) follow from Theorem C.3 (page 36).   ☞

**Example C.5** (power transform metric/snowflake transform metric) [82] Let $(X, \mathsf{d})$ be a *metric space* (Definition 4.5 page 16). $\phi(x) \triangleq x^\alpha, \; \alpha \in (0 : 1]$, is a *metric preserving function* (see Figure 4 page 36 (B))

**Example C.6** ($\alpha$-truncated metric/radar screen metric) [83] Let $(X, \mathsf{d})$ be a *metric space* (Definition 4.5 page 16). $\phi(x) \triangleq \min\{\alpha, x\}, \; \alpha \in \mathbb{R}^+$ is a *metric preserving function* (see Figure 4 page 36 (C)).

**Example C.7** (bounded metric) [84] Let $(X, \mathsf{d})$ be a *metric space* (Definition 4.5 page 16).
$\phi(x) \triangleq \dfrac{x}{1+x}$ is a *metric preserving function* (see Figure 4 page 36 (D)).

**Example C.8** (discrete metric preserving function) [85] Let $\phi$ be a function in $\mathbb{R}^{\mathbb{R}}$.
$\phi(x) \triangleq \left\{ \begin{array}{ll} 0 & \text{for } x \leq 0 \\ 1 & \text{otherwise} \end{array} \right\}$ is a *metric preserving function* (see Figure 4 page 36 (E)).

**Example C.9** Let $\phi$ be a function in $\mathbb{R}^{\mathbb{R}}$.
$\phi(x) \triangleq \left\{ \begin{array}{llll} x & \text{for } 0 \leq x < 1, & 1 & \text{for } 1 \leq x \leq 2, \\ x-1 & \text{for } 2 < x < 3, & 2 & \text{for } x \geq 3 \end{array} \right\}$ is a *metric preserving function* (see Figure 4 page 36 (F)).

---

# Reference Index

# Subject Index

*Telecommunications Engineering Department, National Chiao-Tung University, Hsinchu, Taiwan;* 國立交通大學 *(Gúo Lì Jiāo Tōng Dà Xúe)* 電信工程學系 *(Diàn Xìn Gōng Chéng Xúe Xì)* 新竹，台灣 *(Xīn Zhú, Tái Wān)*

dgreenhoe@gmail.com